\documentclass[a4paper,11pt]{article}
\usepackage{amsmath}
\usepackage{amssymb}
\usepackage{amsthm}
\usepackage{authblk}
\usepackage{enumitem}
\usepackage[margin = 2cm, top = 2.5cm, bottom = 2.2cm]{geometry}
\usepackage[hidelinks]{hyperref}
\usepackage{cite}
\usepackage{subfigure}
\usepackage{tikz}

\newcommand{\boldx}{\mathbf{x}}

\newcommand{\bolda}{\mathbf{a}}
\newcommand{\boldb}{\mathbf{b}}

\newcommand{\xHn}[1]{H^{#1}}
\newcommand{\xHzero}{\xHn{0}}
\newcommand{\xHone}{\xHn{1}}
\newcommand{\xHtwo}{\xHn{2}}
\newcommand{\xLzero}{L^0}
\newcommand{\xLtwo}{L^2}
\newcommand{\xLinfty}{L^{\infty}}
\newcommand{\xCn}[1]{\mathcal{C}^{#1}}
\newcommand{\xCzero}{\xCn{0}}
\newcommand{\xCone}{\xCn{1}}
\newcommand{\xCtwo}{\xCn{2}}
\newcommand{\xR}{\mathbb{R}}
\newcommand{\xN}{\mathbb{N}}
\newcommand{\xP}{\mathbb{P}}
\newcommand{\xker}{\text{ker}}
\newcommand{\xdif}{\mathrm{d}}
\newcommand{\xdrv}[2]{\frac{\xdif {#1}}{\xdif{#2}}}
\newcommand{\red}[1]{\textcolor{red}{#1}}
\newcommand{\blue}[1]{\textcolor{blue}{#1}}

\newtheorem{thrm}{Theorem}
\newtheorem{lmm}{Lemma}
\newtheorem{dfntn}{Definition}
\newtheorem{rmrk}{Remark}
\newtheorem{crllr}{Corollary}

\begin{document}

\title{Virtual Element Method for the Laplace-Beltrami equation on surfaces}
\author[1]{Massimo Frittelli}
\author[1]{Ivonne Sgura}
\affil[1]{\scriptsize{Dipartimento di Matematica e Fisica “E. De Giorgi”, Università del Salento, via per Arnesano, I-73100 Lecce, Italy}}

\date{}

\maketitle

\begin{abstract} 
We present and analyze a Virtual Element Method (VEM) of arbitrary polynomial order $k\in\xN$ for the Laplace-Beltrami equation on a surface in $\xR^3$. The method combines the Surface Finite Element Method (SFEM) [Dziuk, Elliott, \emph{Finite element methods for surface PDEs}, 2013] and the recent VEM [Beirao da Veiga et al, \emph{Basic principles of Virtual Element Methods}, 2013] in order to handle arbitrary polygonal and/or nonconforming meshes. We account for the error arising from the geometry approximation and extend to surfaces the error estimates for the interpolation and projection in the virtual element function space. In the case $k=1$ of linear Virtual Elements, we prove an optimal $\xHone$ error estimate for the numerical method. The presented method has the capability of handling the typically nonconforming meshes that arise when two ore more meshes are pasted along a straight line. Numerical experiments are provided to confirm the convergence result and to show an application of mesh pasting.
\end{abstract}

\section*{MSC Subject Classification}
65N15, 65N30

\section*{Keywords}
Surface PDEs, Laplace-Beltrami equation, Surface Finite Element Method, Virtual Element Method

\section*{Introduction}
The Virtual Element Method (VEM) is a recent extension of the well-known Finite Element Method (FEM) for the numerical approximation of several classes of partial differential equations on planar domains \cite{beirao2013basic, da2013virtual, mora2015virtual, vacca2015virtual, vacca2016virtual, antonietti2016c, benedetto2016globally}. The main features of the method have been introduced in \cite{beirao2013basic, beirao2014hitchhiker}.

The key feature of VEM is that of being a \emph{polygonal} finite element method, i.e. the method handles elements of quite general polygonal shape, rather than just triangular \cite{beirao2013basic}, and nonconforming meshes \cite{beirao2013basic, de2016nonconforming}. The increased mesh generality provides different advantages, we mention some of them. Nonconforming meshes (i) naturally arise when pasting several meshes to obtain a polygonal approximation of the whole domain \cite{benkemoun2012anisotropic, chen2014memory}, as there is no need to match the nodal points in contrast to conforming pasting techniques \cite{kanai1999interactive, sharf2006snappaste} and (ii) allow simple adaptive refinement strategies \cite{cangiani2014adaptive}. Elements of more general shape and arbitrary number of edges allow (i) flexible approximation of the domain and in particular of its boundary \cite{dai2007n} and (ii) the possibility of enforcing higher regularity to the numerical solution \cite{antonietti2016c, da2013arbitrary, brezzi2013virtual}.

The core idea of the VEM is that, given a polynomial order $k\in\xN$ and a polygonal element $K$, the local basis function space on $K$ includes the polynomials of degree $k$ (thus ensuring the optimal degree of accuracy) plus other basis functions that are not known in closed form \cite{beirao2013basic}. The presence of these \emph{virtual} functions motivates the name of the method. However, the knowledge of certain degrees of freedom attached to the basis functions is sufficient to compute the discrete bilinear forms with a degree of accuracy $k$.

The VEM, introduced for the Laplace equation in two dimensions in the recent publication \cite{beirao2013basic}, has been extended to more complicated PDEs, for example a non exhaustive list is: linear elasticity \cite{da2013virtual}, plate bending \cite{brezzi2013virtual}, fracture problems \cite{benedetto2016globally}, eigenvalue problems \cite{mora2015virtual}, Cahn-Hilliard equation \cite{antonietti2016c}, heat \cite{vacca2015virtual} and wave equations \cite{vacca2016virtual}.

The aim of the present work is to extend the VEM to solve \emph{surface PDEs}, i.e. PDEs having a two-dimensional smooth surface in $\xR^3$ as spatial domain. 
Surface PDEs arise in the modelling of several problems such as advection \cite{flyer2007transport}, water waves \cite{flyer2009radial}, phase separation \cite{tang2005phase}, reaction-diffusion systems and pattern formation \cite{bertalmio2001variational, bergdorf2010lagrangian, barreira2011surface, fuselier2013high, frittelli2016lumped}, tumor growth \cite{chaplain2001spatio}, biomembrane modelling \cite{elliott2010modeling}, cell motility \cite{elliott2012modelling}, superconductivity \cite{du2005approximations}, metal dealloying \cite{eilks2008numerical}, image processing \cite{bertalmio2001variational} and surface modelling \cite{xu2006discrete}.
We will focus on the \emph{Laplace-Beltrami} equation, that is the prototypal second order elliptic PDE on smooth surfaces and corresponds to the extension of the Laplace equation to surfaces \cite[chapter 14]{taylor2013partialvol3}.

Among the various discretisation techniques for surface PDEs existing in literature (see for example \cite{macdonald2009implicit, chaplain2001spatio, fuselier2013high, dziuk2013finite, tuncer2015projected}) we consider the Surface Finite Element Method (SFEM) introduced in the seminal paper \cite{dziuk1988finite}. The core idea is to approximate the surface with a polygonal surface made, as in the planar case, of triangular non-overlapping elements whose vertices belong to the surface and to consider a space of piecewise linear functions. The resulting method is exactly similar to the well-known planar FEM, but the convergence estimates must account for the additional error arising from the approximation of the surface, see \cite{dziuk2013finite} for a thorough analysis of the method. In this paper, we define a Virtual Element Method on polygonal surfaces by combining the approaches of VEM and SFEM, the resulting method will be defined as Surface Virtual Element Method (SVEM). Then we prove, under minimal regularity assumptions on the polygonal mesh, some error estimates for the for the approximation of surfaces and for the projection operators and bilinear forms involved in the method. Furthermore, we prove existence and uniqueness of the discrete solution and a first order (and thus optimal) $\xHone$ error estimate. As an application, we show that the method simply handles composite meshes arising from pasting two (or more) meshes along a straight line.

The structure of the paper is as follows.
In Section \ref{sec:preliminaries} we recall some preliminaries on differential operators and function spaces on surfaces.
In Section \ref{sec:laplacebeltrami} we recall the Laplace-Beltrami equation on arbitrary smooth surfaces without boundary in strong and weak forms.
In Section \ref{sec:spacedisc} we introduce a Virtual Element Method for the Laplace-Beltrami equation, defined on general polygonal approximation of surfaces and for any polynomial order $k\in\xN$.
In Section \ref{estimates} we prove error estimates for the discrete bilinear forms and the approximation of geometry.
In Section \ref{erroranalysis} we prove existence, uniqueness and first order $\xHone$ convergence of the numerical solution.
In Section \ref{sec:pasting} we discuss the application of the method to mesh pasting.
In Section \ref{sec:implementation} we face with the issues related to the implementation of the method.
In Section \ref{sec:numericalexamples} we present two numerical examples to (i) test the order of convergence of the method and (ii) show the application of the method to mesh pasting.

\section{Differential operators on surfaces}
\label{sec:preliminaries}
In this section we recall some fundamental notions concerning surface PDEs. If not explicitly stated, definitions and results are taken from \cite{dziuk2013finite}.
\begin{dfntn}[$\xCn{k}$ surface, normal and conormal vectors]
\label{whatisasurface}
Given $k\in\xN$, a set $\Gamma\subset\xR^3$ is said to be a $\xCn{k}$ \emph{surface} if, for every $\boldx_0\in\Gamma$, there exist an open set $U_{\boldx_0}\subset\xR^3$ containing $\boldx_0$ and a function $\phi_{\boldx_0}\in\xCn{k}(U)$ such that
\begin{equation*}
U_{\boldx_0}\cap \Gamma = \{\boldx\in U_{\boldx_0}|\phi_{\boldx_0}(\boldx)=0\}.
\end{equation*}
The vector field
\begin{equation*}
\nu: \Gamma \rightarrow \xR^3,\ \boldx\mapsto \frac{\nabla\phi_{\boldx}(\boldx)}{\|\nabla\phi_{\boldx}(\boldx)\|}
\end{equation*}
is said to be the \emph{unit normal vector}. We denote by $\partial\Gamma$ the one-dimensional \emph{boundary} of $\Gamma$. If $\partial\Gamma$ has a well-defined tangent direction at each point, the vector field $\mu: \partial \Gamma\rightarrow\xR^3$ such that
\begin{itemize}
\item $\mu(\boldx) \perp \nu(\boldx)\quad \forall \boldx\in\partial \Gamma$;
\item $\mu(\boldx)\perp \partial \Gamma\quad \forall \boldx\in\partial \Gamma$;
\item $\mu(\boldx)$ points outward of $\Gamma$,
\end{itemize}
is called the \emph{conormal unit vector}.
\end{dfntn}

\begin{lmm}[Fermi coordinates]
\label{lemma:fermi}
If $\Gamma$ is a $\xCtwo$ surface, there exists an open set $U\subset\xR^3$ such that every $\boldx\in U$ admits a unique decomposition of the form
\begin{equation}
\label{one-to-one}
\boldx = \bolda(\boldx) + d(\boldx)\nu(\bolda(\boldx)),\qquad\ \bolda(\boldx)\in\Gamma,\ d(\boldx)\in\xR.
\end{equation}
The set $U$ is called the \emph{Fermi stripe} of $\Gamma$ and $(\bolda(\boldx),d(\boldx))$ are called the \emph{Fermi coordinates} of $\boldx$.
\end{lmm}

\begin{dfntn}[Tangential gradient, tangential divergence]
If $\Gamma$ is a $\xCone$ surface, $A$ is an open neighborhood of $\Gamma$ and $f\in\xCone(A,\xR)$, the operator
\begin{equation}
\label{definitiontangentialgradient}
\nabla_\Gamma f : S\rightarrow \xR^3,\ \boldx\mapsto \nabla f(\boldx)-(\nabla f(\boldx)\cdot \nu(\boldx))\nu (\boldx) = P(\boldx)\nabla f (\boldx),
\end{equation}
where $P(\boldx)_{ij} = \delta_{ij}-\nu_i(\boldx)\nu_j(\boldx)$, is called the \emph{tangential gradient} of $f$. The components of the tangential gradient, i.e.
\begin{equation*}
\underline{D}_i f : S\rightarrow \xR, \ \boldx\mapsto P_i(\boldx)\nabla f(\boldx),\qquad i\in\{1,2,3\},
\end{equation*}
where $P_i(\boldx)$ is the $i$-th row of $P(\boldx)$, are called the \emph{tangential derivatives} of $f$. Given a vector field $F\in\xCone(A,\xR^3)$, the operator
\begin{equation*}
\nabla_\Gamma\cdot F: S\rightarrow\xR,\ \boldx\mapsto \sum_{i=1}^3\underline{D}_iF_i(\boldx)
\end{equation*}
is called the \emph{tangential divergence} of $F$.
\end{dfntn}

\begin{thrm}
\label{well-posedness-tangential-derivatives}
Given $\Gamma\subset A$ a $\xCone$ surface, if $f$ and $g$ are $\xCone(A,\xR)$ functions such that $f_{|\Gamma}=g_{|\Gamma}$, then
\begin{equation*}
\nabla_\Gamma f(\boldx) = \nabla_\Gamma g(\boldx)\quad \forall \boldx\in \Gamma. 
\end{equation*}
This means that the tangential gradient of a function only depends on its restriction over $\Gamma$.
\end{thrm}
\noindent
Theorem \ref{well-posedness-tangential-derivatives} makes the following definition well-posed.
\begin{dfntn}[$\xCn{k}(\Gamma)$ functions]
If $\Gamma$ is a $\xCone$ surface, a function $f:\Gamma\rightarrow \xR$ is said to be $\xCone(\Gamma)$ if it is differentiable at any point of $\Gamma$ and its tangential derivatives are continuous over $\Gamma$.\\
If $k\geq 2$ and $\Gamma$ is a $\xCn{k}$ surface, a function $f:\Gamma\rightarrow \xR$ is said to be $\xCn{k}(\Gamma)$ if it is $\xCone(\Gamma)$ and its tangential derivatives are $\xCn{{k-1}}(\Gamma)$ functions.
\end{dfntn}
\begin{dfntn}[Laplace-Beltrami operator]
Given a $\xCtwo$ surface $\Gamma$ and $f\in\xCtwo(\Gamma)$, the operator
\begin{equation*}
\Delta_\Gamma f: \Gamma\rightarrow \xR,\ \boldx\mapsto \nabla_\Gamma\cdot\nabla_\Gamma f(\boldx) = \sum_{i=1}^3\underline{D}_i\underline{D}_i f(\boldx)
\end{equation*}
is called the \emph{Laplace-Beltrami} operator of $f$.
\end{dfntn}
\noindent
We now recall the definitions of some remarkable Sobolev spaces on surfaces.
\begin{dfntn}[Sobolev spaces on surfaces]
Given $s\in\xN$, let $\Gamma$ be a $C^s$ surface and let $\xLzero(\Gamma)$ be the set of measurable functions with respect to the bidimensional Hausdorff measure on $\Gamma$. Consider the \emph{Sobolev spaces}
\begin{align*}
&\xHzero(\Gamma) = \xLtwo(\Gamma) = \left\{f\in \xLzero(\Gamma)\ \middle|\ \int_\Gamma f^2\xdif\sigma < +\infty\right\};\\
&\xHn{r}(\Gamma) = \left\{f\in \xLtwo(\Gamma)\ \middle|\ \underline{D}_i f\in \xHn{{r-1}}(\Gamma)\ \forall i\in\{1,2,3\}\right\},\quad \forall\ 1\leq r\leq s;\\
&\xHn{r}_0(\Gamma) = \left\{f\in \xHn{r}(\Gamma) \middle| \int_\Gamma f = 0 \right\},\quad \forall\ 1\leq r\leq s,
\end{align*}
where derivatives are meant in distributional sense\footnote{See \cite[Chapter 4]{taylor2013partialvol1} or \cite[Definition 2.11]{dziuk2013finite} for a precise definition of distributional tangential derivatives.}. These are Hilbert spaces if endowed with the scalar products
\begin{equation*}
<f,g>_{\xHn{r}(\Gamma)} = \int_\Gamma \left(\sum_{|\boldsymbol{\alpha}|\leq r}\underline{D}_{\boldsymbol{\alpha}}f\underline{D}_{\boldsymbol{\alpha}}g\right)\xdif\sigma\quad \forall f,g\in \xHn{r}(\Gamma)\quad \forall\ 0\leq r\leq s.
\end{equation*}
where $\underline{D}_{\boldsymbol{\alpha}}$ is the multi-index notation for partial tangential derivatives.
\end{dfntn}
\noindent
Norms will be denoted by $\|\cdot\|_{\xLtwo(\Gamma)}$, $\|\cdot\|_{\xHn{r}(\Gamma)}$ and seminorms by $|\cdot|_{\xHn{r}(\Gamma)}$.\\
As well as in the planar case, a Poincar\'{e} inequality holds on $\xHone_0(\Gamma)$.
\begin{thrm}[Poincar\'{e}'s inequality on surfaces]
Given a $\xCtwo$ surface $\Gamma$ with a well-define tangent vector field on the boundary $\partial\Gamma$, there exists $C>0$ such that
\begin{equation}
\label{poincareinequality}
\|f\|_{\xLtwo(\Gamma)} \leq C |f|_{\xHone(\Gamma)}\qquad \forall\ f\in \xHone_0(\Gamma).
\end{equation} 
\end{thrm}
\noindent
A basic result in surface calculus is the following
\begin{thrm}[Green's formula on surfaces]
Given a $\xCtwo$ surface $\Gamma$ with a well-defined tangent vector field on the boundary $\partial\Gamma$ and $f,g\in\xCtwo(\Gamma)$, it holds
\begin{equation}
\label{greenformula}
\int_\Gamma f\Delta_\Gamma g\xdif\sigma = -\int_\Gamma \nabla_\Gamma f\cdot\nabla_\Gamma g\xdif\sigma + \int_{\partial \Gamma} f\frac{\partial g}{\partial\mu}\xdif l,
\end{equation}
where $\frac{\partial g}{\partial \mu}(\boldx) = \nabla_\Gamma g(\boldx)\cdot\mu(\boldx)$ is the \emph{conormal derivative} of $g$ on $\partial\Gamma$.
\end{thrm}

\section{The Laplace-Beltrami equation}
\label{sec:laplacebeltrami}
In this section we introduce the Laplace-Beltrami equation on a surface without boundary, that will be the model problem throughout the paper.\\
Let $\Gamma$ be a $\xCn{3}$ surface without boundary and let $f\in \xLtwo(\Gamma)$ such that $\int_\Gamma f = 0$. Consider the Laplace-Beltrami equation on $\Gamma$, given by 
\begin{equation*}
\begin{cases}
-\Delta_\Gamma  u(\boldx) = f(\boldx),\qquad \boldx\in\Gamma,\\
\hspace{3.5mm}\int_\Gamma  u(\boldx)\xdif\boldx = 0,
\end{cases}
\end{equation*}
and its weak formulation
\begin{equation}
\label{weakform0}
\begin{cases}
&u\in \xHone_0(\Gamma)\\
&\int_\Gamma \nabla_\Gamma u \cdot \nabla_\Gamma \phi = \int_\Gamma f\phi\qquad \forall\ \phi\in \xHone(\Gamma).
\end{cases}
\end{equation}
Notice that, from condition $\int_\Gamma f = 0$, the formulation \eqref{weakform0} is equivalent to
\begin{equation}
\label{weakform}
\begin{cases}
&u\in \xHone_0(\Gamma)\\
&\int_\Gamma \nabla_\Gamma u \cdot \nabla_\Gamma \phi = \int_\Gamma f\phi\qquad \boxed{\forall\ \phi\in \xHone_0(\Gamma).}
\end{cases}
\end{equation}
By considering the bilinear form $a(u,v) := \int_\Gamma \nabla_\Gamma u\cdot \nabla_\Gamma v$ for all $u,v\in \xHone(\Gamma)$ and $\langle u,v\rangle_{\xLtwo(\Gamma)}:= \int_\Gamma uv$ for all $u,v\in \xLtwo(\Gamma)$, \eqref{weakform} becomes
\begin{equation}
\label{weakformulation}
\begin{cases}
&u\in \xHone_0(\Gamma)\\
&a(u,\phi) = \langle f,\phi \rangle_{\xLtwo(\Gamma)}\qquad \forall\ \phi\in \xHone_0(\Gamma).
\end{cases}
\end{equation}
Let us justify the above requirements $\int_{\Gamma}u = 0$ and $\int_{\Gamma} f = 0$. Since $\phi \equiv 1$ is allowed as a test function for the weak Laplace-Beltrami equation \eqref{weakform0}, it follows $\int_\Gamma f = 0$ as a compatibility condition. Moreover, if $u$ fulfills $a(u,\phi) = \langle f,\phi \rangle_{\xLtwo(\Gamma)},\ \forall\phi\in \xHone(\Gamma)$ and $c\in\xR$, then $u+c$ fulfills the same equation; condition $\int_\Gamma u=0$ is thus enforced to provide uniqueness of the solution. Existence and uniqueness for problem \eqref{weakformulation} will be proven rigorously in Theorem \ref{abstracttheorem} in Section \ref{erroranalysis}.

\begin{rmrk}[Surfaces with boundary]
\label{rmk:surfaceswithboundary}
The whole analysis carried out in this paper holds unchanged in the presence of a non-empty boundary, $\partial \Gamma \neq \emptyset$, and homogeneous Neumann boundary conditions. In the case of homogeneous Dirichlet boundary conditions, the analysis still holds if $\xHone_0(\Gamma)$ is the space of $\xHone(\Gamma)$ functions that vanish on $\partial \Gamma$ in a weak sense, see \cite[Chapter 4.5]{taylor2013partialvol1}.
\end{rmrk}

\section{Space discretisation by SVEM}
\label{sec:spacedisc}
In this section, we will address space discretisation of \eqref{weakformulation}. After defining the approximation of geometry and the corresponding discrete function spaces, the Surface Virtual Element Method (SVEM) will be introduced.
  
\subsection{Approximation of the surface}
\label{sec:approximationsurface}
In this section we define a polygonal approximation of the surface $\Gamma$ in Definition \ref{whatisasurface} and a \emph{virtual element} space on this polygonal approximation. The method will thus generalise, in the piecewise linear case, the Surface Finite Element Method (SFEM) \cite{dziuk2013finite} and the Virtual Element Method (VEM) \cite{beirao2013basic} at once.\\
Given a $\xCtwo$ surface $\Gamma$ in $\xR^3$, we constuct a piecewise flat approximate surface $\Gamma_h$, defined as
\begin{equation}
\label{polygonal_approximation}
\Gamma_h = \bigcup_{E \in \mathcal{T}_h}E,
\end{equation}
where
\begin{enumerate}
\item \label{itm:firsttriangulation} $\mathcal{T}_h$ is a finite set of non-overlapping simple polygons, i.e. without holes and with non self-intersecting boundary, having diameters less than or equal to $h>0$;
\item \label{itm:secondtriangulation} $\Gamma_h$ is contained in the Fermi stripe $U$ associated to $\Gamma$, see Lemma \ref{lemma:fermi};
\item \label{itm:thirdtriangulation} $\bolda: \Gamma_h\rightarrow \Gamma$ is one-to-one;
\item \label{itm:fourthtriangulation} the vertices of $\Gamma_h$ lie on $\Gamma$.
\end{enumerate}
Following \cite{dziuk2013finite}, we define how to lift functions from the approximate surface $\Gamma_h$ to the continuous one $\Gamma$.
\begin{dfntn}[Lifted functions]
Let $\Gamma$ be a $\xCtwo$ surface and $\Gamma_h$ be as in \eqref{polygonal_approximation}. Given a function $\phi:\Gamma_h\rightarrow\xR$, its \emph{lift} $\phi^\ell:\Gamma\rightarrow\xR$ is defined by $\phi\circ (\bolda_{|\Gamma_h})^{-1}$. Given a function $\psi:\Gamma\rightarrow\xR$, its \emph{unlift} $\psi^{-\ell}:\Gamma_h\rightarrow\xR$ is defined by $\psi\circ \bolda$.\\
This definition is well-posed thanks to assumption (\ref{itm:thirdtriangulation}).
\end{dfntn}
\noindent
Furthermore, the following mesh regularity requirements will be assumed throughout the paper. There exist $\gamma_1,\gamma_2 > 0$ such that, for all $h>0$ and $E\in\mathcal{T}_h$,
\begin{enumerate}[label = (\Alph*)]
\item \label{itm:firstregularity} $E$ is star-shaped with respect to a ball of radius $\rho_E$ such that
\begin{equation*}
\rho_E \geq \gamma_1 h_E;
\end{equation*}
\item \label{itm:secondregularity} for every pair of nodes $P,Q\in E$, the distance $\|P-Q\|$ fulfills
\begin{equation*}
\|P-Q\| \geq \gamma_2 h_E,
\end{equation*}
\end{enumerate}
where $h_E$ is the diameter of $E$.
\begin{rmrk}
This kind of polygonal approximation has two remarkable subcases:
\label{classicaltriangulations}
\begin{enumerate}
\item \label{itm:firstclassical} If each $E\in\mathcal{T}_h$ has three vertices, we obtain the classical triangulations adopted, for instance, in \cite{dziuk2013finite} and \cite{dziuk1988finite}.
\item \label{itm:secondclassical} If $\Gamma$ is a flat surface, we obtain the polygonal meshes considered in \cite{beirao2013basic}.
\end{enumerate}
\end{rmrk}
\noindent
We remark that the considered class of polygonations includes nonconforming meshes. We will show in Section \ref{sec:pasting} that this feature can be exploited in mesh pasting.

\subsection{Discrete function spaces}
Consider a polynomial degree $k\in\xN$ and $E\in\mathcal{T}_h$. Without loss of generality, $E$ may be assumed to lie in the $(x,y)$ plane. Following \cite{beirao2013basic}, the local virtual space of degree $k$ in $E$ is defined by
\begin{equation}
\label{localspace}
V_h^k(E) = \{v_h\in \xHone(E) \mid {v_h}_{|e}\in\xP_k(e)\ \forall e\ \in \text{edges}(E),\ \Delta v_h \in\xP_{k-2}(E)\}.
\end{equation}
The set of \emph{barycentric polynomials} on $E$
\begin{equation*}
\mathcal{M}_{k-2}(E) = \left\lbrace\left(\frac{\boldx-\boldx_E}{h_E}\right)^{\boldsymbol{\alpha}} \middle | |\boldsymbol{\alpha}|\leq k-2\right\rbrace,
\end{equation*}
where $\boldx_E$ and $h_E$ are the barycenter and the diameter of $E$, respectively, is a basis for $\xP_{k-2}(E)$. For every $v_h\in V_h^k(E)$, the following \emph{degrees of freedom} are defined:
\begin{enumerate}
\item the pointwise value of $v_h$ on the $n_E$ vertices of $E$;
\item the pointwise value of $v_h$ on $k-1$ equally spaced points (different from the vertices) of each edge of $E$;
\item the $\frac{k(k-1)}{2}$ moments
\begin{equation*}
\frac{1}{|E|}\int_E v_h m_{k-2}\qquad \forall\ m_{k-2}\in\mathcal{M}_{k-2}(E).
\end{equation*}
\end{enumerate}
In \cite{beirao2013basic} it has been proven that these degrees of freedom are unisolvent for the space $V_h^k(E)$ in \eqref{localspace}. Given $w\in \xHn{s}(E)$, $s>1$, by Sobolev's embedding theorem we have $w\in \xCzero(E)$. Hence, the degrees of freedom of $w$ are well-defined. If $n^{dof}_E$ is the number of degrees of freedom on $E$ and $dof_i(v_h)$ denotes the $i^{th}$ degree of freedom of $v_h$, $i=1,\dots,n^{dof}_E$, the unique function $w_I\in V_h^k(E)$ such that
\begin{equation}
\label{interpolationdefinition}
dof_i(w-w_I)=0\qquad \forall i=1,\dots, n^{dof}_E
\end{equation}
is said to be the interpolant of $w$. If $\tilde{E} = \bolda(E)$ is the curved triangle corresponding to $E$ and $v\in \xHtwo(\tilde{E})$, we have $v^{-\ell}\in \xHtwo(E)$, see Theorem \ref{equivalence}. The unique function $v_I\in V_h^k(E)$ such that
\begin{equation}
\label{liftedinterpolationdefinition}
dof_i(v^{-\ell}-v_I)=0\qquad \forall i=1,\dots, n^{dof}_E
\end{equation}
is said to be the interpolant of $v$.
\begin{rmrk}
From the definition of $V_h(E)$ in \eqref{localspace} it follows that
\begin{enumerate}
\item $\xP_k(E) \subseteq V_h^k(E)$;
\item every $v_h\in V_h^k(E)$ is explicitly known on the boundary $\partial E$, but not on the interior $\overset{\circ}{E}$;
\item if $k=1$, $V_h^1(E)$ is the set of harmonic functions in $E$ being piecewise linear on the boundary $\partial E$ and the set of local degrees of freedom collapses to the pointwise values on the vertices of $E$;
\item if $k=1$ and $E$ is a triangle, $V_h^1(E) = \xP_1(E)$; this is the only case in which $V_h^k(E) = \xP_k(E)$, i.e the VEM method reduces to FEM. 
\end{enumerate}
\end{rmrk}
\noindent
The global discrete space will be defined by
\begin{equation*}
V_h^k = \{v_h\in\xCzero(\Gamma_h) \mid {v_h}_{|E} \in V_h^k(E)\ \forall\ E\in\mathcal{T}_h\}.
\end{equation*}
Furthermore, we define the zero-averaged virtual space $W_k^k$ by
\begin{equation}
\label{zeroaveragedspace}
W_h^k = \left\{v_h\in V_h^k \ \middle|\ \int_{\Gamma_h} v_h=0\right\}.
\end{equation}
Finally, we define the following broken $H^s$ seminorms, $s\in \{1,2\}$, on the polygonal surface $\Gamma_h$:
\begin{equation*}
|v_h|_{h,s} = \sqrt{\sum_{E\in\mathcal{T}_h}|v_{h|E}|_{\xHn{s}(E)}^2}\qquad \forall v_h\in \prod_{E\in\mathcal{T}_h}\xHn{s}(E)
\end{equation*}

\subsection{The Surface Virtual Element Method}
\label{SVEM}
We may write a discrete formulation for \eqref{weakformulation}:
\begin{equation}
\label{noncomputableformulation}
\begin{cases}
&u_h\in W_h^k\\
&\int_{\Gamma_h}\nabla_{\Gamma_h} u_h\cdot\nabla_{\Gamma_h}\phi_h = \int_{\Gamma_h}f_I\phi_h \qquad \forall\ \phi_h\in W_h^k,
\end{cases}
\end{equation}
where $f_I$ is the interpolant of $f$ defined piecewise in \eqref{liftedinterpolationdefinition}.
\begin{rmrk}[Regularity of $f$]
In the following we assume $f\in\xHtwo(\Gamma)$, such that, from Sobolev's embedding theorem, its pointwise values (and thus its interpolant $f_I$) are well-defined. We remark that, in the framework of surface PDEs, the problem of numerically handling $\xHn{s}(\Gamma)$, $0\leq s\leq 1$, load terms is intrinsically challenging. In fact, if the pointwise values of $f$ are not available, then any approximation $\bar{f}$ of $f$ defined on $\Gamma_h$ must account for the mapping $\bolda:\Gamma_h\rightarrow \Gamma$ in \eqref{one-to-one} that, in general, is not computable.
\end{rmrk}
\noindent
By introducing the bilinear forms $\bar{a}(u_h,v_h) := \int_{\Gamma_h}\nabla_{\Gamma_h} u_h\cdot\nabla_{\Gamma_h}v_h$ for all $u_h,v_h\in V_h^k$ and $\langle U_h,V_h\rangle_{\xLtwo(\Gamma_h)}$ for all $U_h, V_h\in \xLtwo(\Gamma_h)$, problem \eqref{noncomputableformulation} is equivalent to
\begin{equation}
\label{discreteformulation}
\begin{cases}
&u_h\in W_h^k\\
&\bar{a}(u_h,\phi_h) = \langle f_I,\phi_h\rangle_{\xLtwo(\Gamma_h)} \qquad \forall\ \phi_h\in W_h^k.
\end{cases}
\end{equation}
We recall that functions in $V_h^k$ are \emph{virtual}, i.e. they are not known explicitly, then $\bar{a}(\cdot,\cdot)$ and $\langle f_I,\cdot \rangle_{\xLtwo(\Gamma_h)}$ are not computable. We thus need to write a computable approximation of problem \eqref{discreteformulation}. To this end, following \cite{beirao2013basic}, an approximate bilinear form $a_h(\cdot,\cdot)$ and an approximate linear form $\langle f_h, \cdot\rangle_h$ will be constructed instead of $\bar{a}(\cdot,\cdot)$ and $\langle f_I,\cdot\rangle_{\xLtwo(\Gamma_h)}$, respectively.\\
Given the following decomposition of $\bar{a}$
\begin{equation*}
\bar{a}(v_h,w_h) = \sum_{E\in\mathcal{T}_h} \bar{a}_E({v_h}_{|E}, {w_h}_{|E})\qquad \forall\ v_h,w_h\in V_h^k,
\end{equation*}
consider the projection $\Pi^\nabla_E: V_h^k(E)\rightarrow \xP_k(E)$ defined by
\begin{equation}
\label{pnabla}
\begin{cases}
&\Pi^\nabla_E(v_h)\in\xP_k(E)\\
&\bar{a}_E\left(\Pi^\nabla_E(v_h), q_k\right) = \bar{a}_E(v_h, q_k)\qquad \forall\ q_k\in\xP_k(E)
\end{cases}
\end{equation}
together with
\begin{equation}
\label{additionalcondition}
\begin{cases}
&\displaystyle \sum_{P\in\ \text{nodes}(E)} \Pi^\nabla_E v_h(P) = \sum_{P\in\ \text{nodes}(E)} v_h(P)\qquad \text{if}\ k=1,\\
&\displaystyle \int_E\Pi^\nabla_E v_h = \int_E v_h\qquad \text{if}\ k>1.
\end{cases}
\end{equation}
The additional condition \eqref{additionalcondition} is enforced to fix the free constant in \eqref{pnabla}. We now prove that $\Pi^\nabla_E$ is computable. The right hand side of \eqref{pnabla} is computable, since
\begin{equation}
\label{computability}
\bar{a}_E(v_h, w_h) = \int_E \nabla v_h\cdot \nabla q_k = -\underbrace{\int_E v_h\Delta q_k}_{\text{(term 1)}} + \underbrace{\int_{\partial E} v_h(\nabla q_k\cdot \nu_E)}_{\text{(term 2)}},
\end{equation}
where $\nu_E$ is the unit outward vector on $\partial E$, and
\begin{itemize}
\item $\Delta q_{k}$ is a polynomial of degree $k-2$, thus (term 1) in \eqref{computability} is a linear combination of the moments of $v_h$;
\item $v_h$ and $q_k$ are both explicitly known (and polynomials) on $\partial E$.
\end{itemize}
The right hand side of \eqref{additionalcondition} is computable since
\begin{itemize}
\item for $k=1$, it is a combination of the pointwise values, $v_h(P)$, $P\in\text{nodes}(E)$, that are degrees of freedom;
\item for $k>1$, it is one of the moments of $v_h$, as the space of barycentric monomials $\mathcal{M}_{k-2}$ contains the constant monomial $m_{k-2}\equiv 1$.
\end{itemize}
Hence, $\Pi^\nabla_E$ is computable. By expressing $v_h,w_h\in V_h^k(E)$ as
\begin{align}
&v_h = \Pi^\nabla_E v_h + (I-\Pi^\nabla_E)v_h;\\
&w_h = \Pi^\nabla_E w_h + (I-\Pi^\nabla_E)v_w,
\end{align}
the form $\bar{a}_E(v_h,w_h)$ may be decomposed as
\begin{equation}
\label{decomposition}
\bar{a}_E(v_h,w_h) = \underbrace{\bar{a}_E(\Pi^\nabla_E v_h,\Pi^\nabla_E w_h)}_{\text{(term 1)}} + \underbrace{\bar{a}_E((I-\Pi^\nabla_E)v_h, (I-\Pi^\nabla_E)w_h)}_{\text{(term 2)}} \qquad \forall\ v_h,w_h\in V_h^k(E),
\end{equation}
because cross-terms vanish due to the definition of $\Pi^\nabla_E$. Notice that (term 1) in \eqref{decomposition} is computable since $\Pi^\nabla_E$ is computable, but does not scale as $\bar{a}_E$ on $\xker(\Pi^\nabla_E)$. Then (term 2) in \eqref{decomposition} cannot be neglected, but must be approximated in a suitable way. To this end we recall that, under the regularity assumptions \ref{itm:firstregularity}-\ref{itm:secondregularity}, the bilinear form
\begin{equation}
\label{stabilisingform}
S_E(v_h,w_h) = \sum_{i=1}^{n^{dof}_E} dof_i(v_h)dof_i(w_h)\qquad \forall\ v_h,w_h\in V_h^k(E),
\end{equation}
scales as $\bar{a}_E$ on the kernel of $\Pi^\nabla_E$, i.e. there exist $c^{*} > c_{*} > 0$ such that
\begin{equation}
\label{stabilisingform_stability}
c_{*}\bar{a}_E(v_h,v_h) \leq S_E(v_h,v_h) \leq c^{*}\bar{a}_E(v_h,v_h)\qquad \forall\ v_h\in \xker(\Pi^\nabla_E),
\end{equation}
see \cite{beirao2013basic}. Consider now a local approximate form $a_{h,E}$ defined by
\begin{equation}
\label{localapproximateform}
a_{h,E}(v_h,w_h) =\ \bar{a}_E(\Pi^\nabla_E v_h,\Pi^\nabla_E w_h) + S_E((I-\Pi^\nabla_E)v_h, (I-\Pi^\nabla_E)w_h)\qquad\forall\ v_h,w_h\in V_h^k(E).
\end{equation}
Notice that, since $\Pi^\nabla_E q_k = q_k$ for all $q_k\in \xP_k(E)$, the local form \eqref{localapproximateform} satisfes the consistency property
\begin{equation}
\label{consistency}
a_{h,E}(v_h,q_k) = \bar{a}_E(v_h,q_k)\qquad \forall v_h\in V_h^k(E),\qquad \forall q_k\in\xP_k(E).
\end{equation}
A global approximate gradient form defined by pasting the local ones:
\begin{equation}
\label{globalapproximateform}
a_h(v_h,w_h) = \sum_{E\in\mathcal{T}_h} a_{h,E}({v_h}_{|E}, {w_h}_{|E})\qquad\forall\ v_h,w_h\in V_h^k.
\end{equation}
We want to define an approximate $\xLtwo$ form and the approximate right hand side. For $n\geq 0$ and for every $E\in\mathcal{T}_h$, consider the local $\xLtwo(E)$ projection $\Pi^E_{n}: V_h^k(E)\rightarrow \xP_{n}(E)$ given by
\begin{equation}
\label{L2projection}
\begin{cases}
\Pi^E_{n}(v_h)\in \xP_{n}(E),\\
\langle \Pi^E_{n} v_h, q_{n} \rangle_{\xLtwo(E)} = \langle v_h, q_{n} \rangle_{\xLtwo(E)}\qquad \forall q_{n}\in \xP_{n}(E).
\end{cases}
\end{equation}
In \eqref{L2projection} we choose $n$ depending on $k$ as
\begin{equation*}
n_k = \begin{cases}
0\hspace*{1.33cm} \text{if}\qquad k=1;\\
k-2\qquad \text{if}\qquad k\geq 2.
\end{cases}
\end{equation*}
We remark that
\begin{itemize}
\item for $k=1$, we have that $\Pi^E_0 = \Pi^\nabla_E$ (see for instance \cite{beirao2014hitchhiker}), hence $\Pi_0^E$ is computable;
\item for $k\geq 2$, $\Pi^E_{k-2}$ is computable since, in \eqref{L2projection}, $\langle v_h, q_{k-2} \rangle_{\xLtwo(E)}$ is a linear combination of the moments of $v_h$.
\end{itemize}
Following \cite{beirao2014hitchhiker}, and in analogy with the approximate gradient form \eqref{localapproximateform}, we consider the following local approximate $\xLtwo$ form:
\begin{equation*}
\langle v_h,w_h\rangle_{\xLtwo_{h,E}} := \int_{\Gamma_h} \Pi^E_{n_k} v_h\Pi^E_{n_k} w_h + S_E((I-\Pi^E_{n_k})v_h, (I-\Pi^E_{n_k})w_h),\qquad \forall v_h, w_h\in V_h^k(E),
\end{equation*}
where $S_E$ and $\Pi^E_{n_k}$ are defined in \eqref{stabilisingform} and \eqref{L2projection}, respectively. Notice that the approximate $\xLtwo$ form \eqref{approximateL2form} fulfills the consistency property.
\begin{equation*}
\langle v_h,q_{n_k}\rangle_{\xLtwo_{h,E}} = \langle v_h,q_{n_k} \rangle_{\xLtwo(E)},\qquad \forall v_h\in V_h^k(E),\qquad \forall q_{n_k}\in\xP_{n_k}(E). 
\end{equation*}
As a consequence, for any $k\in\mathbb{N}$ we have that
\begin{equation}
\label{exactintegral}
\langle v_h,1\rangle_{\xLtwo_h, E} = \int_{E} v_h\qquad \forall v_h\in V_h^k(E),
\end{equation}
i.e. the integral of any $V_h^k$ function can be computed exactly. A computable global approximate $\xLtwo$ form is obtained by pasting the local ones:
\begin{equation}
\label{approximateL2form}
\langle v_h,w_h\rangle_{\xLtwo_h} = \sum_{E\in\mathcal{T}_h} \langle {v_h}_{|E}, {w_h}_{|E}\rangle_{\xLtwo_{h,E}},\qquad \forall v_h,w_h\in V_h^k.
\end{equation}
Property \eqref{exactintegral} implies that the space $W_h^k$ defined in \eqref{zeroaveragedspace} can be represented as
\begin{equation}
\label{zeroaveragedrepresentation}
W_h^k = \{v_h\in V_h | \langle v_h,1\rangle_{\xLtwo_h} = 0\},
\end{equation}
hence $W_k$ is computable. To approximate the right hand side, following \cite{beirao2013basic}, for any function $g\in \xHone(\Gamma_h)$ we consider the functional $\langle g,\cdot\rangle_h$ defined by
\begin{equation}
\label{computablefunctional}
\begin{split}
\langle g,v_h\rangle_h =
\begin{cases}
\displaystyle\sum_{E\in\mathcal{T}_h}\int_E g \sum_{V\in\text{nodes}(E)}\dfrac{v_h(V)}{n_E}\qquad \text{if}\ k=1\\
\displaystyle\sum_{E\in\mathcal{T}_h}\int_{E} \Pi_{k-2}^E g\ v_h \hspace*{2.1cm} \text{if}\ k\geq 2
\end{cases}\qquad 
\forall\ v_h\in V_h^k,
\end{split}
\end{equation}
where $n_E$ is the number of nodes of $E$. From \eqref{exactintegral} we have that, if $g\in V_h^k$, $\langle g, v_h\rangle_h$ is computable, given the degrees of freedom of $g$ and $v_h$. Furthermore, notice that
\begin{equation}
\label{exactintegralh}
\langle g, 1\rangle_h = \int_{\Gamma_h} g,\qquad \forall g\in \xHone(\Gamma).
\end{equation}
In order to simplify the  implementation, as we will discuss in Section \ref{sec:implementation}, we define an approximate computable load term $f_h := f_I - \frac{\langle f_I,1\rangle_{\xLtwo_h}}{|\Gamma_h|}$. From \eqref{exactintegral} it follows that $f_h$ is zero averaged and, from \eqref{exactintegralh}, $f_h$ fulfills
\begin{equation}
\label{discretezeroaveragedloadterm}
\langle f_h,1\rangle_h = 0
\end{equation}
We may now write a computable discrete problem:
\begin{equation}
\label{computablediscreteformulation}
\begin{cases}
&u_h\in W_{h}^k\\
&a_h(u_h,\phi_h) = \langle f_h,\phi_h\rangle_h \qquad \forall\ \phi_h\in W_{h}^k.
\end{cases}
\end{equation}
We will discuss the implementation of \eqref{computablediscreteformulation} in Section \ref{sec:implementation}.\\
The error analysis will be carried out in the following steps:
\begin{enumerate}
\item the geometric and interpolation error estimates in \cite{dziuk2013finite} will be extended to our polygonal/virtual setting;
\item the error between the continuous weak formulation \eqref{weakformulation} and the computable discrete one \eqref{computablediscreteformulation} will be estimated by extending the analogous convergence theorem in \cite{beirao2013basic}.
\end{enumerate}
In Section \ref{estimates} we deal with step (1), in Section \ref{erroranalysis} we deal with step (2).

\section{Interpolation, projection and geometric error estimates}
\label{estimates}
We start this section by recalling some results from \cite{beirao2013basic}. The following theorem addresses the projection error on $\xP_k(E)$, $E\in\mathcal{T}_h$.
\begin{thrm}
\label{vem_interp}
Under the regularity assumption \ref{itm:firstregularity}, there exists $C>0$, depending only on $k$ and $\gamma$, such that for every $1\leq s\leq k+1$ and for all $w\in \xHn{s}(E)$ there exists a $w_\pi\in\xP_k(E)$ such that
\begin{equation}
\label{projectionest}
\|w-w_\pi\|_{\xLtwo(E)} + h_E|w-w_\pi|_{\xHone(E)} \leq C h_E^s|w|_{\xHn{s}(E)}
\end{equation}
\qed
\end{thrm}
\noindent
We now address interpolation in $V_k(E)$, $E\in\mathcal{T}_h$. The following theorem from \cite{beirao2013basic} gives an interpolation error estimate in $V_h^k(E)$.
\begin{thrm}
\label{interpolationbasic}
Under the regularity assumption \ref{itm:firstregularity}, there exists $C>0$, depending only on $k$ and $\gamma$, such that for every $2 \leq s \leq k+1$ and for all $w\in \xHn{s}(E)$, the interpolant $w_I\in V_h^k(E)$ satisfies
\begin{equation}
\label{interpolationestimate}
\|w-w_I\|_{\xLtwo(E)} + h|w-w_I|_{\xHone(E)} \leq Ch_E^s|w|_{\xHn{s}(E)}.
\end{equation}
\qed
\end{thrm}
\noindent
To approximate a function $u:\Gamma\rightarrow\xR$ with an interpolant defined on $\Gamma_h$ as in \eqref{liftedinterpolationdefinition}, also a geometric error has to be taken into account. The following lemma generalizes Lemma 4.1 in \cite{dziuk2013finite} to our polygonal approximations of the surface.
\begin{lmm}
\label{geometricestimates}
Let $\Gamma_h$ be a polygonal approximation of $\Gamma$ as in \eqref{polygonal_approximation}. For any $E\in\mathcal{T}_h$, the oriented distance function introduced in \eqref{one-to-one} fulfills
\begin{equation}
\label{distanceestimate}
\|d\|_{\xLinfty(E)} \leq Ch^2.
\end{equation}
The quotient $\delta_h$ between the smooth and the discrete surface measures defined by $\xdif A = \delta_h\xdif A_h$ satisfies
\begin{equation}
\label{quotientestimate}
\|1-\delta_h\|_{\xLinfty(\Gamma_h)} \leq Ch^2.
\end{equation}
Let $P$ and $P_h$ be the projections onto the tangent planes of the smooth and the discrete surfaces, respectively, that is $P_{ij} = \delta_{ij}-\nu_i\nu_j$, $P_{h,ij} = \delta_{ij}-\nu_{h,i}\nu_{h,j}$, and define
\begin{equation}
\label{rhestimate}
R_h = \frac{1}{\delta_h}P(I-d\mathcal{H})P_h(I-d\mathcal{H}),
\end{equation}
where $\mathcal{H}$ is the Weingarten map defined by $\mathcal{H}_{ij} = \xdrv{\nu_i}{x_j}$. Then
\begin{equation}
\label{finalgeomestimate}
\|(I-R_h)P\|_{\xLinfty(\Gamma_h)} \leq Ch^2.
\end{equation}
In all of the claimed inequalities $C$ depends only on the curvature of $\Gamma$.
\begin{proof}
Consider $E\in\mathcal{T}_h$, see Fig. \ref{fig:proof_step0}. First of all we prove that
\begin{equation}
\label{firsttoprove}
\|d\|_{\xLinfty(\partial E)} \leq \frac{h^2}{8}|d|_{\xCtwo(E)}.
\end{equation}
To this end, let $\boldx_B\in\partial E$ and let $e$ be an edge of $E$ such that $\boldx_B\in e$, see Fig. \ref{fig:proof_step1}. Then, if $d_e$ is the linear interpolant of $d$ on $e$, from the classical Lagrange interpolation estimates and from the fact that $d_e \equiv 0$ since $d$ vanishes at the endpoints of $e$, we have
\begin{equation*}
|d(\boldx_B)| \leq \|d_e\|_{\xLinfty(e)} + \|d-d_e\|_{\xLinfty(e)} \leq \frac{|e|^2}{8}|d|_{\xCtwo(e)} \leq \frac{h^2}{8}|d|_{\xCtwo(E)},
\end{equation*}
that proves \eqref{firsttoprove}. Now let $\boldx\in \overset{\circ}{E}$ and let $s$ be any straight line contained in the plane of $E$ and passing through $\boldx$, let $\boldx_1,\boldx_2 \in s\cap\partial E$ such that $[\boldx_1,\boldx_2]\subset E$, see Fig. \ref{fig:proof_step2}, and let $d_s$ be the linear interpolant of $d$ on $s$. From classical Lagrange interpolation estimates and from \eqref{firsttoprove} we have that
\begin{equation}
\label{secondtoprove}
\begin{split}
|d(\boldx)| &\leq \|d_s\|_{\xLinfty(s)} + \|d-d_s\|_{\xLinfty(s)} = \max (|d(\boldx_1)|,|d(\boldx_2)|) + \|d-d_s\|_{\xLinfty(s)}\\
 &\leq \|d\|_{\xLinfty(\partial E)} + \frac{|s|^2}{8}|d|_{\xCtwo(s)} \leq \frac{h^2}{8}|d|_{\xCtwo(E)} + \frac{h^2}{8}|d|_{\xCtwo(E)} \leq \frac{h^2}{4}|d|_{\xCtwo(E)} \leq \frac{h^2}{4}|d|_{\xCtwo(U)},
\end{split}
\end{equation}
where $U$ is the Fermi stripe of $\Gamma$, but $|d|_{\xCtwo(U)}$ depends only on the curvature of $\Gamma$, thus \eqref{secondtoprove} proves \eqref{distanceestimate}. To prove \eqref{quotientestimate}, \eqref{rhestimate} and \eqref{finalgeomestimate}, we proceed as in Lemma 4.1 in \cite{dziuk2013finite} using estimate \eqref{distanceestimate} for polygonal meshes.
\end{proof}
\end{lmm}
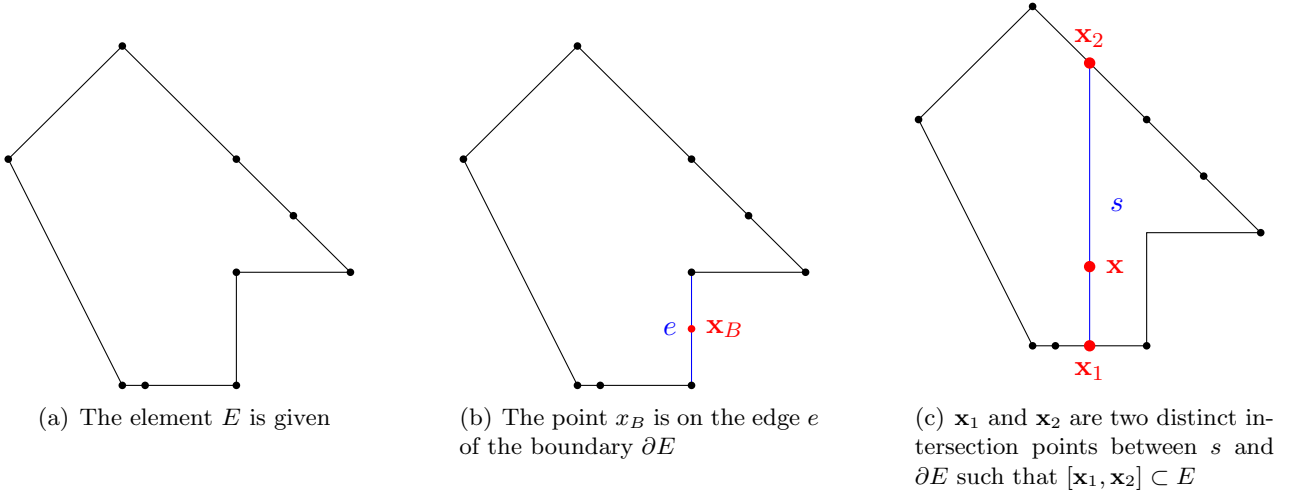
\begin{figure}
\begin{center}
\subfigure[The element $E$ is given]{\begin{tikzpicture}[scale=1.5]

\draw (0,0)--(1,0)--(1,1)--(2,1)--(0,3)--(-1,2)--cycle;

\node[shape=circle, fill = black, minimum size = 1mm, inner sep = 0pt] at (0,0) {};
\node[shape=circle, fill = black, minimum size = 1mm, inner sep = 0pt] at (1,0) {};
\node[shape=circle, fill = black, minimum size = 1mm, inner sep = 0pt] at (2,1) {};
\node[shape=circle, fill = black, minimum size = 1mm, inner sep = 0pt] at (-1,2) {};
\node[shape=circle, fill = black, minimum size = 1mm, inner sep = 0pt] at (0,3) {};
\node[shape=circle, fill = black, minimum size = 1mm, inner sep = 0pt] at (1,1) {};

%HANGING NODES
\node[shape=circle, fill = black, minimum size = 1mm, inner sep = 0pt] at (1,2) {};
\node[shape=circle, fill = black, minimum size = 1mm, inner sep = 0pt] at (1.5,1.5) {};
\node[shape=circle, fill = black, minimum size = 1mm, inner sep = 0pt] at (0.2,0) {};

\end{tikzpicture}\label{fig:proof_step0}}
\hspace{1cm}
\subfigure[The point $x_B$ is on the edge $e$ of the boundary $\partial E$]{\begin{tikzpicture}[scale=1.5]

\draw (1,1)--(2,1)--(0,3)--(-1,2)--(0,0)--(1,0);
\draw[color = blue]  (1,0)--(1,1);

\node[shape=circle, fill = black, minimum size = 1mm, inner sep = 0pt] at (0,0) {};
\node[shape=circle, fill = black, minimum size = 1mm, inner sep = 0pt] at (1,0) {};
\node[shape=circle, fill = black, minimum size = 1mm, inner sep = 0pt] at (1,1) {};
\node[shape=circle, fill = black, minimum size = 1mm, inner sep = 0pt] at (2,1) {};
\node[shape=circle, fill = black, minimum size = 1mm, inner sep = 0pt] at (-1,2) {};
\node[shape=circle, fill = black, minimum size = 1mm, inner sep = 0pt] at (0,3) {};

%HANGING NODES
\node[shape=circle, fill = black, minimum size = 1mm, inner sep = 0pt] at (1,2) {};
\node[shape=circle, fill = black, minimum size = 1mm, inner sep = 0pt] at (1.5,1.5) {};
\node[shape=circle, fill = black, minimum size = 1mm, inner sep = 0pt] at (0.2,0) {};

\node[shape=circle, fill = red, minimum size = 1mm, inner sep = 0pt, label=right:{\red{$\boldx_B$}}, label = left:{\blue{$e$}}] at (1,0.5) {};
\end{tikzpicture}\label{fig:proof_step1}}
\hspace{1cm}
\subfigure[$\boldx_1$ and $\boldx_2$ are two distinct intersection points between $s$ and $\partial E$ such that {$[\boldx_1,\boldx_2 ] \subset E$}]{\begin{tikzpicture}[scale=1.5]

\draw (0,0)--(1,0)--(1,1)--(2,1)--(0,3)--(-1,2)--cycle;
\draw[color = blue]  (0.5,0)--(0.5,2.5);

\node[label = right:{\blue{$s$}}] at (0.5,1.25) {};

\node[shape=circle, fill = black, minimum size = 1mm, inner sep = 0pt] at (0,0) {};
\node[shape=circle, fill = black, minimum size = 1mm, inner sep = 0pt] at (1,0) {};
\node[shape=circle, fill = black, minimum size = 1mm, inner sep = 0pt] at (2,1) {};
\node[shape=circle, fill = black, minimum size = 1mm, inner sep = 0pt] at (-1,2) {};
\node[shape=circle, fill = black, minimum size = 1mm, inner sep = 0pt] at (0,3) {};
\node[shape=circle, fill = red, minimum size = 1.5mm, inner sep = 0pt, label = right:{\red{$\boldx$}}] at (0.5,0.7) {};
\node[shape=circle, fill = red, minimum size = 1.5mm, inner sep = 0pt, label = below:{\red{$\boldx_1 $}}] at (0.5,0) {};
\node[shape=circle, fill = red, minimum size = 1.5mm, inner sep = 0pt, label = above:{\red{$\boldx_2 $}}] at (0.5,2.5) {};

%HANGING NODES
\node[shape=circle, fill = black, minimum size = 1mm, inner sep = 0pt] at (1,2) {};
\node[shape=circle, fill = black, minimum size = 1mm, inner sep = 0pt] at (1.5,1.5) {};
\node[shape=circle, fill = black, minimum size = 1mm, inner sep = 0pt] at (0.2,0) {};

\end{tikzpicture}\label{fig:proof_step2}}
\end{center}
\caption{Some steps of the proof of Lemma \ref{geometricestimates}.}\label{fig:euclideangeometry}
\end{figure}
\noindent
The following lemma generalizes Lemma 4.2 in \cite{dziuk2013finite} to our polygonal setting and provides lower and upper bounds for some norms of arbitrary functions when they are unlifted from $\Gamma$ to $\Gamma_h$ or lifted from $\Gamma_h$ to $\Gamma$. In particular, it provides an equivalence between the $\xLtwo(\Gamma_h)$ and $\xLtwo(\Gamma)$ norms and between the $\xHone(\Gamma_h)$ and $\xHone(\Gamma)$ seminorms.

\begin{lmm}
\label{equivalence}
Let $w:\Gamma_h\rightarrow\xR$ with lift $w^\ell:\Gamma\rightarrow\xR$. Let $\bolda: \Gamma_h\rightarrow\Gamma$ be the projection onto $\Gamma$ defined in \eqref{one-to-one} and, for every $E\in\mathcal{T}_h$, let $\tilde{E} = \bolda (E)\subset\Gamma$ be the curved triangle corresponding to $E\in\mathcal{T}_h$. Then
\begin{align}
\label{equivalence1}
&\frac{1}{C} \|w\|_{\xLtwo(E)} \leq \|w^\ell\|_{\xLtwo(\tilde{E})} \leq C\|w\|_{\xLtwo(E)};\\
\label{equivalence2}
&\frac{1}{C} \|\nabla_E w\|_{\xLtwo(E)} \leq \|\nabla_{\tilde{E}}w^\ell\|_{\xLtwo(\tilde{E})} \leq C\|\nabla_E w\|_{\xLtwo(E)};\\
\label{equivalence3}
&\|\nabla^2_E w\|_{\xLtwo(E)} \leq C\|\nabla^2_{\tilde{E}}w^\ell\|_{\xLtwo(\tilde{E})} + Ch_E\|\nabla_{\tilde{E}}w^\ell\|_{\xLtwo(\tilde{E})},
\end{align}
if the norms exist, where $C$ depends only on the surface area and the curvature of $\Gamma$.
\begin{proof}
We use the estimates of Lemma \ref{geometricestimates} and proceed exactly as in \cite[Lemma 4.2]{dziuk2013finite}.
\end{proof}
\end{lmm}
\noindent
The following result provides, in the case $k=1$, error estimates for the interpolation in $(V_h^1)^\ell$ and the projection on $\left(\prod_{E}\xP_1(E)\right)^\ell$. The interpolation result extends to SVEM Lemma 4.3 in \cite{dziuk2013finite} for the triangular SFEM.
\begin{thrm}
\label{interpolationlemma}
Given a $\xCtwo$ surface $\Gamma$, there exists $C>0$ such that, for all $v\in \xHtwo(\Gamma)$ and $w\in \xHn{s}(\Gamma)$, $s\in\{1,2\}$, and for all $h>0$, then
\begin{itemize}
\item the interpolant $v_I\in V_h^1$ fulfills
\begin{equation}
\label{interpolationlemmaestimate}
\|v-v_I^\ell\|_{\xLtwo(\Gamma)} + h|v-v_I^\ell |_{\xHone(\Gamma)} \leq Ch^2\left(|v|_{\xHtwo(\Gamma)} + h|v|_{\xHone(\Gamma)}\right);
\end{equation}
\item there exists a projection $w_\pi\in \prod_{E}\xP_1(E)$ such that
\begin{equation}
\label{projectionlemmaestimate}
\|w-w_\pi^\ell\|_{\xLtwo(\Gamma)} + h|w-w_I^\ell |_{h,1} \leq Ch^s\left(|w|_{\xHn{s}(\Gamma)} + h|w|_{\xHone(\Gamma)}\right).
\end{equation}
\end{itemize}
\begin{proof}
From Lemma \ref{equivalence}, $w^{-\ell}\in \xHone(\Gamma_h)\cap\prod_{E} \xHn{s}(E)$. Let $w_\pi$ be the $\prod_{E}\xP_1(E)$ projection of $w^{-\ell}$ as in \eqref{projectionest} and let $v_I$ be the $V_h^1$ interpolant of $v^{-\ell}$ defined piecewise by \eqref{interpolationdefinition}. From Theorems \ref{vem_interp} and \ref{interpolationbasic}, by summing piecewise contributions, we have
\begin{align}
\label{projecstep}
&\|w^{-\ell}-w_\pi\|_{\xLtwo(\Gamma_h)} + h|w^{-\ell}-w_\pi|_{h,1} \leq Ch^s|w^{-\ell}|_{2,h},\\
\label{interpstep}
&\|v^{-\ell}-v_I\|_{\xLtwo(\Gamma_h)} + h|v^{-\ell}-v_I|_{\xHone(\Gamma_h)} \leq Ch^2 |v^{-\ell}|_{2,h}.
\end{align}
From \eqref{projecstep}, \eqref{interpstep} and Lemma \ref{equivalence} we have
\begin{align}
&\|w-w_\pi^\ell\|_{\xLtwo(\Gamma)} + h|w-w_\pi^\ell|_{h,1} \leq Ch^s\left(|w|_{\xHn{s}(\Gamma)} + h|w|_{\xHone(\Gamma)}\right),\\
&\|v-v_I^\ell\|_{\xLtwo(\Gamma)} + h|v-v_I^\ell|_{\xHone(\Gamma)} \leq Ch^2\left(|v|_{\xHtwo(\Gamma)} + h|v|_{\xHone(\Gamma)}\right),
\end{align}
that are the desired estimates.
\end{proof}
\end{thrm}

\begin{rmrk}
\label{interpolationremark}
Let $\Gamma$ be a sufficiently smooth surface. If the equivalences \eqref{equivalence1}-\eqref{equivalence3} are guaranteed for arbitrary $\xHn{k}$ norms, the estimate \eqref{interpolationlemmaestimate} immediately generalizes to
\begin{equation*}
|w-w_I^\ell|_{\xHn{r}(\Gamma)} \leq C h^{k+1-r}|w|_{\xHn{{k+1}}(\Gamma)},\qquad r=0,\dots,k,
\end{equation*}
where $k$ is the order of the VEM space defined in Equation \eqref{localspace}.
\end{rmrk}
\noindent
The following Lemma generalizes Lemma 4.7 in \cite{dziuk2013finite} to our polygonal/virtual setting and provides bounds for the geometric errors in the bilinear forms.
\begin{lmm}
\label{bilinearforms}
For any $(v,w)\in \xHone(\Gamma_h)\times \xHone(\Gamma_h)$, the following estimates hold:
\begin{align}
\label{bilinearform1}
|\langle v^\ell,w^\ell\rangle_{\xLtwo(\Gamma)} - \langle v,w\rangle_{\xLtwo(\Gamma_h)}| \leq Ch^2\|v^\ell\|_{\xLtwo(\Gamma)}\|w^\ell\|_{\xLtwo(\Gamma)};\\
\label{bilinearform2}
|a(v^\ell,w^\ell)-\bar{a}(v,w)| \leq Ch^2\|\nabla_\Gamma v^\ell\|_{\xLtwo(\Gamma)}\|\nabla_\Gamma w^\ell\|_{\xLtwo(\Gamma)},
\end{align}
where $C$ depends only on the geometry of $\Gamma$.
\begin{proof}
We proceed as in Lemma 4.7 of \cite{dziuk2013finite}, but here using the generalized estimates \eqref{distanceestimate}-\eqref{finalgeomestimate} given in the previous Lemma \ref{geometricestimates}.
\end{proof}
\end{lmm}
\noindent
In the first section we have recalled the Poincar\'{e} inequality \eqref{poincareinequality} in $\xHone_0(\Gamma)$. In the following theorem we prove an analogous inequality in $\xHone_0(\Gamma_h)$, i.e on polygonal surfaces $\Gamma_h$ of the type \eqref{polygonal_approximation}.
\begin{thrm}[Poincar\'{e} inequality in $\xHone_0(\Gamma_h)$]
Let $\Gamma$ be a closed $\xCtwo$ orientable surface in $\xR^3$. Then there exist $h_0>0$ and $C>0$ depending on $\Gamma$ such that, for all $0<h<h_0$ and $\Gamma_h$ as in \eqref{polygonal_approximation},
\begin{equation}
\label{triangulatedpoincare}
\|v\|_{\xLtwo(\Gamma_h)} \leq C|v|_{\xHone(\Gamma_h)}\qquad \forall v\in \xHone_0(\Gamma_h).
\end{equation}
\begin{proof}
From \eqref{equivalence1} and the triangle inequality we have
\begin{equation}
\label{poincare0}
\|v\|_{\xLtwo(\Gamma_h)} \leq C \|v^\ell\|_{\xLtwo(\Gamma)} \leq C\left(\left\|v^\ell-\frac{1}{|\Gamma|}\int_\Gamma v^\ell\right\|_{\xLtwo(\Gamma)} + \frac{1}{|\Gamma|^\frac{1}{2}}\int_\Gamma v^\ell\right).
\end{equation}
Now, from \eqref{equivalence1} we have that $v^\ell-\frac{1}{|\Gamma|}\int_\Gamma v^\ell \in \xHone_0(\Gamma)$. Then, from Poincar\'{e}'s inequality \eqref{poincareinequality} and \eqref{equivalence2} it follows that
\begin{equation}
\label{poincare1}
\left\|v^\ell-\frac{1}{|\Gamma|}\int_\Gamma v^\ell\right\|_{\xLtwo(\Gamma)} \leq C |v^\ell|_{\xHone(\Gamma)} \leq C |v|_{\xHone(\Gamma_h)}.
\end{equation}
Furthermore, from \eqref{equivalence1}, \eqref{bilinearform1} and the fact that $v$ is zero-averaged on $\Gamma_h$, it follows that
\begin{equation}
\label{poincare2}
\frac{1}{|\Gamma|^\frac{1}{2}}\int_\Gamma v^\ell  \leq  \frac{1}{|\Gamma|^\frac{1}{2}}\left( \left|\int_{\Gamma_h} v\right| + Ch^2\|v^\ell\|_{\xLtwo(\Gamma)}|\Gamma|^{\frac{1}{2}}\right) \leq Ch^2\|v\|_{\xLtwo(\Gamma_h)}.
\end{equation}
Combining \eqref{poincare0}, \eqref{poincare1} and \eqref{poincare2} we have
\begin{equation*}
(1-Ch^2)\|v\|_{\xLtwo(\Gamma_h)} \leq C|v|_{\xHone(\Gamma)}.
\end{equation*}
By choosing, for instance, $h_0 = \frac{1}{\sqrt{2C}}$, the result follows.
\end{proof}
\end{thrm}
\noindent
Concerning the convergence rates of the above results we observe that:
\begin{itemize}
\item As shown in Lemma \ref{bilinearforms}, in the approximation of the bilinear forms \eqref{bilinearform1} and \eqref{bilinearform2}, the polygonal approximation of geometry yields a geometric error that is quadratic in $\xLtwo$ norm and linear in $\xHone$ norm. In fact, this Lemma is based on the geometric estimates of Lemma \ref{geometricestimates}.
\item The interpolation error on $\Gamma$, as shown by \eqref{interpolationlemmaestimate} in Lemma \ref{interpolationlemma} (and its proof) arises from two sources. The first one is the interpolation error on flat polygons (cp. Lemma \ref{vem_interp}). The second one is given by the geometric estimates given in Lemma \ref{geometricestimates}. Nevertheless, the order of accuracy only depends on the first of these two sources, and thus on the VEM order $k$ (cp. Remark \ref{interpolationremark}).
\end{itemize}
This rate gap implies that choosing $k>1$ in \eqref{localspace} will not improve the convergence rate of the method, since geometric error dominates over the interpolation one. For this reason, in what follows, we restrict our study to the case $k=1$, i.e. ``piecewise linear virtual elements''. The same drawback occurs with the standard SFEM \cite{dziuk2013finite} of higher order, $k>1$; in \cite{demlow2009higher} it has been shown that a finite element space of degree $k$ defined on a suitable curvilinear triangulation of degree $k$ (isoparametric elements) provides a SFEM with the same convergence rate as polynomial interpolation of degree $k$. This suggests that, to formulate a SVEM of order $k>1$, a different approximation of the surface is needed. We close this section proving an error estimate for the approximate right hand side $<f_h,v_h>$ in the discrete formulation \eqref{computablediscreteformulation} for $k=1$.
\begin{thrm}
Let $f\in \xHone_0(\Gamma)$. Under the regularity assumptions \ref{itm:firstregularity}-\ref{itm:secondregularity}, there exists $C>0$ depending on $\Gamma$ and $\gamma$ such that
\begin{equation}
\label{righthanderror}
|\langle f,v_h^\ell\rangle_{\xLtwo(\Gamma)} - \langle f_h, v_h\rangle_h | \leq Ch\left(|f|_{\xHone(\Gamma)} + h|f|_{\xHtwo(\Gamma)}\right)|v_h^\ell|_{\xHone(\Gamma)}\qquad \forall v_h\in W_h^1.
\end{equation}
\begin{proof}
Let $f_I$ be as in \eqref{noncomputableformulation} and $f_h$ be as in \eqref{computablediscreteformulation}. We split the error as
\begin{equation}
\label{splittederror}
\begin{split}
|\langle f,v_h^\ell\rangle_{\xLtwo(\Gamma)} - \langle f_h, v_h\rangle_h | &\leq |\langle f,v_h^\ell\rangle_{\xLtwo(\Gamma)} - \langle f_I, v_h\rangle_{\xLtwo(\Gamma_h)}| + |\langle f_I, v_h \rangle_{\xLtwo(\Gamma_h)} - \langle f_h, v_h\rangle_{\xLtwo(\Gamma_h)}|\\
&+ |\langle f_h, v_h \rangle_{\xLtwo(\Gamma_h)} - \langle f_h, v_h\rangle_h|.
\end{split}
\end{equation}
From the Cauchy-Schwarz inequality,\eqref{bilinearform1} we obtain
\begin{equation}
\begin{split}
|\langle f,v_h^\ell\rangle_{\xLtwo(\Gamma)} - \langle f_I, v_h\rangle_{\xLtwo(\Gamma_h)}| &\leq |\langle f-f_I^\ell,v_h^\ell\rangle_{\xLtwo(\Gamma)}| + |\langle f_I^\ell,v_h^\ell\rangle_{\xLtwo(\Gamma)} - \langle f_I, v_h\rangle_{\xLtwo(\Gamma_h)}|\\
&\leq \|f-f_I^\ell\|_{\xLtwo(\Gamma)}\|v_h^\ell\|_{\xLtwo(\Gamma)} + Ch^2\|f_I^\ell\|_{\xLtwo(\Gamma)}\|v_h^\ell\|_{\xLtwo(\Gamma)}
\end{split}
\end{equation}
From the Cauchy-Schwarz inequality, the definition of $f_h$ and \eqref{bilinearform1} we have
\begin{equation}
\begin{split}
&|\langle f_I, v_h \rangle_{\xLtwo(\Gamma_h)} - \langle f_h, v_h\rangle_{\xLtwo(\Gamma_h)}| \leq |\Gamma_h|^{-\frac{1}{2}}|\langle f_I,1\rangle_{\xLtwo(\Gamma_h)}| \|v_h\|_{\xLtwo(\Gamma_h)}\\
\leq &|\Gamma_h|^{-\frac{1}{2}}\left(|\langle f_I^\ell-f,1\rangle_{\xLtwo(\Gamma)}| +Ch^2\|f_I^\ell\|_{\xLtwo(\Gamma)}\right)\|v_h\|_{\xLtwo(\Gamma_h)}\\
\leq &\left(\|f_I^\ell-f\|_{\xLtwo(\Gamma)} + Ch^2\|f_I^\ell\|_{\xLtwo(\Gamma)}\right)\|v_h\|_{\xLtwo(\Gamma_h)}
\end{split}
\end{equation}
Following \cite{beirao2013basic}, we know that
\begin{equation}
|\langle f_h, v_h \rangle_{\xLtwo(\Gamma_h)} - \langle f_h, v_h\rangle_h| \leq Ch|f_h|_{1,h}|v_h|_{\xHone(\Gamma_h)},
\end{equation}
but, from the definition of $f_h$ and from \eqref{equivalence2} it follows that
\begin{equation}
\label{splitted1}
|f_h|_{1,h} = |f_I|_{\xHone(\Gamma_h)} \leq C|f_I^\ell|_{\xHone(\Gamma)}.
\end{equation}
Combining \eqref{splittederror}-\eqref{splitted1}, using \eqref{equivalence1},\eqref{equivalence2}, \eqref{interpolationlemmaestimate}, the Poincar\'{e} inequalities \eqref{poincareinequality}, \eqref{triangulatedpoincare} and the triangle inequality we obtain
\begin{equation*}
\begin{split}
&|\langle f,v_h^\ell\rangle_{\xLtwo(\Gamma)} - \langle f_h, v_h\rangle_h | \leq \left(\|f-f_I^\ell\|_{\xLtwo(\Gamma)} + Ch|f_I^\ell|_{\xHone(\Gamma)} + Ch^2\|f_I^\ell\|_{\xLtwo(\Gamma)}\right)|v_h^\ell|_{\xHone(\Gamma)}\\
\leq &\left((1+Ch^2)\|f-f_I^\ell\|_{\xLtwo(\Gamma)} + Ch^2\|f\|_{\xLtwo(\Gamma)} + Ch|f-f_I^\ell|_{\xHone(\Gamma)} + Ch|f|_{\xHone(\Gamma)} \right)|v_h^\ell|_{\xHone(\Gamma)}\\
\leq &C\left((h^2+h^4)|f|_{\xHtwo(\Gamma)} + (h+h^3+h^5)|f|_{\xHone(\Gamma)} \right)|v_h^\ell|_{\xHone(\Gamma)} \leq Ch\left(|f|_{\xHone(\Gamma)} + h|f|_{\xHtwo(\Gamma)}\right)|v_h^\ell|_{\xHone(\Gamma)},
\end{split}
\end{equation*}
that is the desired estimate.
\end{proof}
\end{thrm}

\section{Existence, uniqueness and error analysis}
\label{erroranalysis}
The following theorem, that is the main result of this paper, extends Theorem 3.1 in \cite{beirao2013basic} for the VEM on planar domains to the Laplace-Beltrami equation on surfaces. In fact, it provides: (i) the existence and the uniqueness of the solution for both the continuous \eqref{weakformulation} and the discrete problem \eqref{computablediscreteformulation} and (ii) an abstract convergence result. As a corollary, an optimal $\xHone(\Gamma)$ error estimate for problem \eqref{computablediscreteformulation} will be given.
\begin{thrm}[Abstract convergence theorem]
\label{abstracttheorem}
Let $a:\xHone_0(\Gamma)\times \xHone_0(\Gamma)\rightarrow\xR$ be the bilinear form defined by
\begin{equation*}
a(u,v) = \int_\Gamma \nabla_\Gamma u\cdot\nabla_\Gamma v\qquad \forall\ u,v\in \xHone_0(\Gamma),
\end{equation*}
and let $a_h: W_h^1\times W_h^1\rightarrow\xR$ be any symmetric bilinear form such that 
\begin{equation}
\label{discretedecomposition}
a_h(u_h,v_h) = \sum_{E\in\mathcal{T}_h}a_{h,E}(u_{h|E},v_{h|E})
\end{equation}
where, for all $E\in\mathcal{T}_h$, $a_{h,E}$ is a symmetric bilinear form on $V_h^1(E)\times V_h^1(E)$ such that
\begin{align}
\label{quasiconsistency}
|&a_{h,E}(p,v_{h,E}) - a_{\tilde{E}}(p^\ell,v_{h,E}^\ell)| \leq Ch^2|p^\ell|_{\xHone(\tilde{E})}|v_{h,E}^\ell|_{\xHone(\tilde{E})} \qquad \forall\ v_{h,E}\in V_h^1(E)\ \forall\ p\in\xP_1(E);\\
\label{bilinearequivalence}
&\alpha_{*}a_{\tilde{E}}(v_{h,E}^\ell,v_{h,E}^\ell) \leq a_{h,E}(v_{h,E},v_{h,E}) \leq \alpha^{*}a_{\tilde{E}}(v_{h,E}^\ell,v_{h,E}^\ell) \qquad \forall\ v_{h,E}\in V_h^1(E),
\end{align}
where $\alpha_{*}$ and $\alpha^{*}$ are independent of $h$ and $E\in\mathcal{T}_h$.\\
Let $F\in \xLtwo(\Gamma)'$ and $F_h\in (W_h^1)'$ be linear continuous functionals. Consider the problems
\begin{align}
\label{theorem_continuousproblem}
&\begin{cases}
&{u}\in \xHone_0(\Gamma)\\
&a({u},v) = {F}(v)\quad \forall v\in \xHone_0(\Gamma)
\end{cases}\\
\label{theorem_discreteproblem}
&\begin{cases}
&u_h\in W_h\\
&a_h(u_h,v_h) = F_h(v_h)\quad \forall v_h\in W_h^1
\end{cases}
\end{align}
Both of these problems have a unique solution and the following error estimate holds
\begin{equation}
\label{errorestimate}
|u-u_h^\ell|_{\xHone(\Gamma)} \leq C \left(|u-u_\pi^\ell|_{h,1} + |u-u_I^\ell|_{\xHone(\Gamma)} + \mathcal{F}_h + h\|F\|_{\xLtwo(\Gamma)'}\right),
\end{equation}
where $\mathcal{F}_h$ is the smallest constant such that
\begin{equation}
\label{definitionfh}
|F(v_h^\ell)-F_h(v_h)| \leq \mathcal{F}_h|v_h^\ell|_{\xHone(\Gamma)}\qquad \forall v_h\in W_h^1.
\end{equation}
\begin{proof}
Existence and uniqueness follow from Lax-Milgram's theorem. In fact, from the Poincar\'{e} inequality \eqref{poincareinequality} on $\xHone_0(\Gamma)$, the bilinear form $a$ is coercive and, from the Cauchy-Schwarz inequality, it is continuous. The bilinear form $a_h$ is coercive since
\begin{equation*}
\begin{split}
|a_h(v_h,v_h)| &\underset{\eqref{discretedecomposition}}{=} \left| \sum_{E\in\mathcal{T}_h} a_{h,E}(v_{h|E},v_{h|E})\right| \underset{\eqref{bilinearequivalence}}{\geq} \alpha_* \sum_{E\in\mathcal{T}_h} \left| a_{\tilde{E}}(v_{h|E}^\ell,v_{h|E}^\ell)\right| = \alpha_*  \sum_{E\in\mathcal{T}_h} |v_h^\ell|_{\xHone(\tilde{E})}^2\\
&\underset{\eqref{equivalence2}}{\geq} C \sum_{E\in\mathcal{T}_h} |v_h|_{\xHone(E)}^2 = C|v_h|_{\xHone(\Gamma_h)}^2 \underset{\eqref{triangulatedpoincare}}{\geq} C \|v_h\|_{\xHone(\Gamma_h)}^2,
\end{split}
\end{equation*}
for all $v_h\in W_h^1$. Now we prove that $a_h$ is continuous. To this end, since $a_h$ is symmetric and coercive (i.e. positive definite), then it fulfills the Cauchy-Schwarz inequality. Then we have
\begin{equation*}
\begin{split}
& |a_h(v_h,w_h)| \underset{\eqref{discretedecomposition}}{\leq} \sum_{E\in\mathcal{T}_h} |a_{h,E}(v_{h|E},w_{h|E})|\leq \sum_{E\in\mathcal{T}_h} a_{h,E}(v_{h|E},v_{h|E})^\frac{1}{2}a_{h,E}(w_{h|E},w_{h|E})^\frac{1}{2}\\
\underset{\eqref{bilinearequivalence}}{\leq} &\alpha^* \sum_{E\in\mathcal{T}_h} a_{\tilde{E}}(v_{h|E}^\ell,v_{h|E}^\ell)^\frac{1}{2}a_{\tilde{E}}(w_{h|E}^\ell,w_{h|E}^\ell)^\frac{1}{2} = \alpha^* \sum_{E\in\mathcal{T}_h} |v_{h}^\ell|_{\xHone(\tilde{E})} |w_{h}^\ell|_{\xHone(\tilde{E})}\\
\underset{\eqref{equivalence2}}{\leq} &C  \sum_{E\in\mathcal{T}_h}  |v_{h}|_{\xHone(E)} |w_{h}|_{\xHone(E)} \leq C \left(\sum_{E\in\mathcal{T}_h} |v_h|_{\xHone(E)}^2\right)^\frac{1}{2}\left(\sum_{E\in\mathcal{T}_h} |w_h|^2_{\xHone(E)}\right)^\frac{1}{2}\\
= & C |v_h|_{\xHone(\Gamma_h)}|w_h|_{\xHone(\Gamma_h)} \leq C \|v_h\|_{\xHone(\Gamma_h)}\|w_h\|_{\xHone(\Gamma_h)},
\end{split}
\end{equation*}
for all $v_h\in W_h^1$. Hence, problems \eqref{theorem_continuousproblem} and \eqref{theorem_discreteproblem} meet the assumptions of Lax-Milgram's theorem.\\
If $\bolda: \Gamma_h\rightarrow\Gamma$ is the projection onto $\Gamma$ defined in \eqref{one-to-one}, then for any $E\in\mathcal{T}_h$, let $\tilde{E} = \bolda(E)$ be the curved triangle corresponding to $E$. Let $u_\pi\in \prod_{E\in\mathcal{T}_h}\xP_1(E)$ be the projection of $u$ defined in \eqref{projectionlemmaestimate} and let $u_I\in W_h^1$ be the interpolant of $u$ defined in \eqref{interpolationlemmaestimate}. From \cite[Theorem 3.3]{dziuk2013finite}, The solution of \eqref{theorem_continuousproblem} fulfills $u\in \xHtwo(\Gamma)$ and thus $u_\pi$ and $u_I$ are well-defined. Let $\delta_h = u_h-u_I$. It holds that
\begin{equation*}
\begin{split}
&\alpha_{*}|\delta_h^\ell|_W^2 = \alpha_{*}a(\delta_h^\ell,\delta_h^\ell) \leq a_h(\delta_h,\delta_h) =a_h(u_h,\delta_h) - a_h(u_I,\delta_h)\\ &\underset{\eqref{discretedecomposition}}{=} F_h(\delta_h) - \sum_{E\in\mathcal{T}_h} a_{h,E}(u_I,\delta_h) =F_h(\delta_h) - \sum_{E\in\mathcal{T}_h}\left(a_{h,E}(u_I-u_\pi,\delta_h)+a_{h,E}(u_\pi,\delta_h)\right)\\
&\underset{\eqref{quasiconsistency}}{\leq} F_h(\delta_h) - \sum_{E\in\mathcal{T}_h}\left(a_{h,E}(u_I-u_\pi,\delta_h)+a_{\tilde{E}}(u_\pi^\ell,\delta_h^\ell)\right) + Ch^2\sum_{E\in\mathcal{T}_h}|u_\pi^\ell|_{\xHone(\tilde{E})}|\delta_h^\ell|_{\xHone(\tilde{E})}\\
&=F_h(\delta_h) - \sum_{E\in\mathcal{T}_h}\left(a_{h,E}(u_I-u_\pi,\delta_h)+a_{\tilde{E}}(u_\pi^\ell-u,\delta_h^\ell) +a_{\tilde{E}}(u,\delta_h^\ell)\right) + Ch^2\left(|u_\pi^\ell|_{h,1}^2 + |\delta_h^\ell|_{\xHone(\Gamma)}^2\right)\\
&= F_h(\delta_h)-a(u,\delta_h^\ell) - \sum_{E\in\mathcal{T}_h}\left( a_{h,E}(u_I-u_\pi,\delta_h) + a_{\tilde{E}} (u_\pi^\ell-u,\delta_h^\ell)\right) + Ch^2\left(|u_\pi^\ell|_{h,1}^2 + |\delta_h^\ell|_{\xHone(\Gamma)}^2\right)\\
&=F_h(\delta_h)-F(\delta_h^\ell) - \sum_{E\in\mathcal{T}_h}\left(a_{h,E}(u_I-u_\pi,\delta_h) + a_{\tilde{E}}(u_\pi^\ell-u,\delta_h^\ell)\right) + Ch^2\left(|u_\pi^\ell|_{h,1}^2 + |\delta_h^\ell|_{\xHone(\Gamma)}^2\right).
\end{split}
\end{equation*}
From \eqref{bilinearequivalence}, \eqref{definitionfh} and the continuity of $a$ and $a_h$ we obtain
\begin{equation*}
(\alpha_{*}-Ch^2)|\delta_h^\ell|_{\xHone(\Gamma)}^2 \leq \mathcal{F}_h|\delta_h^\ell|_{\xHone(\Gamma)} + |u_I-u_\pi|_{h,1}|\delta_h|_{\xHone(\Gamma_h)} + |u_\pi^\ell-u|_{h,1}|\delta_h^\ell|_{\xHone(\Gamma)} + Ch^2|u_\pi^\ell|_{h,1}^2
\end{equation*}
For $h$ sufficiently small, by exploiting \eqref{equivalence2}, we obtain
\begin{equation}
\label{algebraicinequality}
|\delta_h^\ell|_{\xHone(\Gamma)}^2 \leq C(\mathcal{F}_h + |u_I^\ell-u_\pi^\ell|_{h,1} + |u_\pi^\ell-u|_{h,1})|\delta_h^\ell|_{\xHone(\Gamma)} + Ch^2|u_\pi^\ell|_{h,1}^2
\end{equation}
By defining $A = \mathcal{F}_h + |u_I^\ell-u_\pi^\ell|_{h,1} + |u_\pi^\ell-u|_{h,1}$ and solving the second-degree-algebraic inequality \eqref{algebraicinequality} we have
\begin{equation*}
|\delta_h^\ell|_{\xHone(\Gamma)} \leq \frac{CA}{2} + \frac{1}{2}\sqrt{C^2A^2+4Ch^2|u_\pi^\ell|_{h,1}^2} \leq \frac{CA}{2} + \frac{1}{2}(CA + 2\sqrt{C}h|u_\pi^\ell|_{h,1}) \leq CA + Ch|u_\pi^\ell|_{h,1}.
\end{equation*}
By recalling the definition of $A$ and applying the triangle inequality, we get
\begin{equation*}
|u-u_h^\ell|_{\xHone(\Gamma)} \leq C(\mathcal{F}_h + |u-u_I^\ell|_{\xHone(\Gamma)} + |u-u_\pi^\ell|_{h,1}) + Ch|u_\pi^\ell|_{h,1}.
\end{equation*}
By applying the triangle inequality to the last term, we obtain
\begin{equation*}
|u-u_h^\ell|_{\xHone(\Gamma)} \leq C\left(\mathcal{F}_h + |u-u_I^\ell|_{\xHone(\Gamma)} + (1+h)|u-u_\pi^\ell|_{h,1} + h|u|_{\xHone(\Gamma)}\right).
\end{equation*}
The obvious stability estimate $|u|_{\xHone(\Gamma)}\leq C\|F\|_{\xLtwo(\Gamma)'}$, where $C$ is the constant in the Poincar\'{e} inequality \eqref{poincareinequality}, together with $h\leq h_0$, complete the proof.
\end{proof}
\end{thrm}
\noindent
From the abstract framework given in Theorem \ref{abstracttheorem} we are now ready to derive the $\xHone(\Gamma)$ error estimate between the continuous problem \eqref{weakformulation} and the discrete one \eqref{computablediscreteformulation}.
\begin{crllr}[$\xHone(\Gamma)$ error estimate]
\label{crllr:errorestimate}
Problems \eqref{weakformulation} and \eqref{computablediscreteformulation} have a unique solution. Let $u$ and $u_h$ be the their solutions, respectively. Under the mesh regularity assumptions \ref{itm:firstregularity}-\ref{itm:secondregularity}, if $f\in \xHtwo_0(\Gamma)$, the following estimate holds:
\begin{equation}
\label{convergenceestimate}
|u-u_h^\ell|_{\xHone(\Gamma)} \leq Ch(|u|_{\xHtwo(\Gamma)} + |f|_{\xHone(\Gamma)}) + Ch^2 |f|_{\xHtwo(\Gamma)}.
\end{equation}
\begin{proof}
In Theorem \ref{abstracttheorem}, we choose
\begin{align*}
&F(v) = \langle f,v\rangle_{\xLtwo(\Gamma)} \qquad\ \forall v\in \xHone(\Gamma);\\
&F_h(v_h) = \langle f_h,v_h\rangle_h \qquad \ \forall v_h\in W_h^1,
\end{align*}
with $a_h$ defined in \eqref{localapproximateform}, \eqref{globalapproximateform}. Under the regularity assumptions \ref{itm:firstregularity}-\ref{itm:secondregularity},
\begin{enumerate}
\item Assumption \eqref{quasiconsistency} follows from \eqref{consistency} and \eqref{bilinearform2};
\item Assumption \eqref{bilinearequivalence} follows from \eqref{decomposition}, \eqref{stabilisingform_stability}, \eqref{localapproximateform} and \eqref{equivalence2};
\item From \cite[Theorem 3.3]{dziuk2013finite} we have $u\in \xHtwo(\Gamma)$. Then, Theorem \ref{interpolationlemma} provides
\begin{equation}
\label{errorbound1}
|u-u_\pi^\ell|_{h,1} + |u-u_I^\ell|_{\xHone(\Gamma)} < Ch (|u|_{\xHtwo(\Gamma)}+h|u|_{\xHone(\Gamma)});
\end{equation}
\item if $f \in \xHone_0(\Gamma)$, The Poincar\'{e} inequality \eqref{poincareinequality} provides
\begin{equation}
\label{errorbound2}
\|F\|_{\xLtwo(\Gamma)'} = \|f\|_{\xLtwo(\Gamma)} \leq C|f|_{\xHone(\Gamma)},
\end{equation}
and \eqref{righthanderror} yields
\begin{equation}
\label{errorbound3}
\mathcal{F}_h\leq Ch(|f|_{\xHone(\Gamma)}+h|f|_{\xHtwo(\Gamma)})
\end{equation}
\end{enumerate}
By plugging \eqref{errorbound1}-\eqref{errorbound3}, into the abstract error bound \eqref{errorestimate}, we obtain
\begin{equation}
\label{errorbound4}
|u-u_h^\ell|_{\xHone(\Gamma)} \leq Ch(|u|_{\xHtwo(\Gamma)} + |f|_{\xHone(\Gamma)}) + Ch^2(|u|_{\xHone(\Gamma)} + |f|_{\xHtwo(\Gamma)})
\end{equation}
By plugging the Poincar\'{e} inequality \eqref{poincareinequality}, the stability estimate $|u|_{\xHone(\Gamma)} \leq C\|F\|_{\xLtwo(\Gamma)'}$ and \eqref{errorbound2} into \eqref{errorbound4}, the result follows.
\end{proof}
\end{crllr}

\section{Pasting polygonal surfaces along a straight line}
\label{sec:pasting}
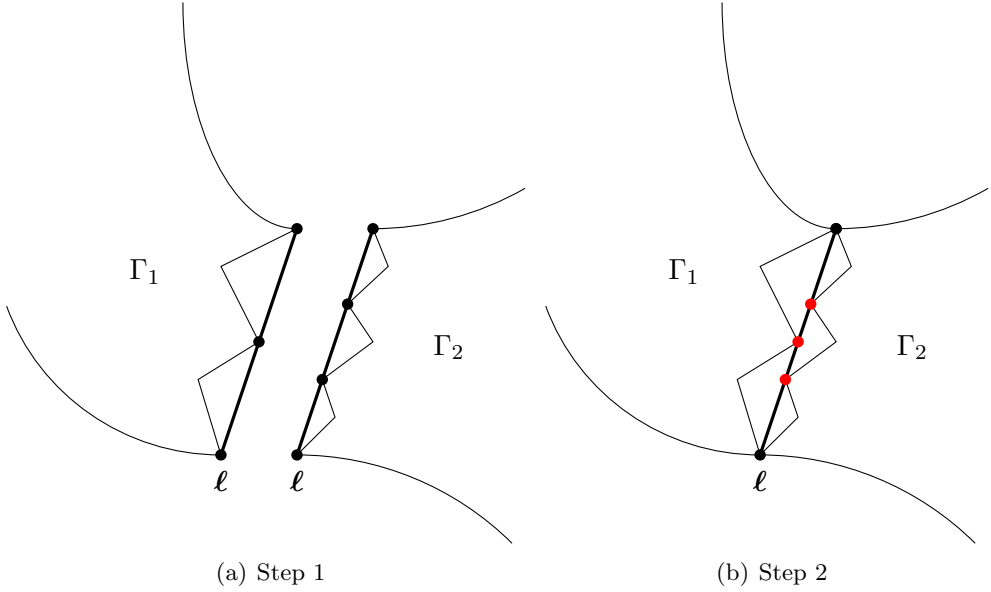
\begin{figure}
\begin{center}
\subfigure[Step 1]{\begin{tikzpicture}
\draw[very thick] (0,0)--(1,3);
\draw[very thick] (1,0)--(2,3);

\draw (1,3) arc (270:180:1.5 and 3);
\draw (0,0) arc (270:200:3);
\draw (2,3) arc (270:300:4);
\draw (1,0) arc (90:45:4);

\node[label = {$\Gamma_1$}] at (-1,2) {};
\node[label = {$\Gamma_2$}] at (3,1) {};
\node[shape=circle, fill = black, minimum size = 1.5mm, inner sep = 0pt, label = below:{$\boldsymbol{\ell}$}] at (0,0) {};
\node[shape=circle, fill = black, minimum size = 1.5mm, inner sep = 0pt, label = below:{$\boldsymbol{\ell}$}] at (1,0) {};

\draw (0,0)--(0.5,1.5)--(-0.3,1)--(0,0);
\draw (0.5,1.5)--(0,2.5)--(1,3);
\node[shape=circle, fill = black, minimum size = 1.5mm, inner sep = 0pt] at (0.5,1.5) {};
%\node[shape=circle, fill = black, minimum size = 1.5mm, inner sep = 0pt] at (-0.3,1) {};
%\node[shape=circle, fill = black, minimum size = 1.5mm, inner sep = 0pt] at (0,2.5) {};
\node[shape=circle, fill = black, minimum size = 1.5mm, inner sep = 0pt] at (1,3) {};

\draw (1,0)--(1.333,1)--(1.5,0.5)--(1,0);
\draw (1.333,1)--(1.666,2)--(2,1.5)--(1.333,1);
\draw (1.666,2)--(2,3)--(2.2,2.5)--(1.666,2);
\node[shape=circle, fill = black, minimum size = 1.5mm, inner sep = 0pt] at (1.333,1) {};
%\node[shape=circle, fill = black, minimum size = 1.5mm, inner sep = 0pt] at (1.5,0.5) {};
\node[shape=circle, fill = black, minimum size = 1.5mm, inner sep = 0pt] at (1.666,2) {};
%\node[shape=circle, fill = black, minimum size = 1.5mm, inner sep = 0pt] at (2,1.5) {};
%\node[shape=circle, fill = black, minimum size = 1.5mm, inner sep = 0pt] at (2.2,2.5) {};
\node[shape=circle, fill = black, minimum size = 1.5mm, inner sep = 0pt] at (2,3) {};

\end{tikzpicture}}
\subfigure[Step 2]{\begin{tikzpicture}
\draw[very thick] (0,0)--(1,3);

\draw (1,3) arc (270:180:1.5 and 3);
\draw (0,0) arc (270:200:3);
\draw (1,3) arc (270:300:4);
\draw (0,0) arc (90:45:4);

\node[label = {$\Gamma_1$}] at (-1,2) {};
\node[label = {$\Gamma_2$}] at (2,1) {};
\node[shape=circle, fill = black, minimum size = 1.5mm, inner sep = 0pt, label = below:{$\boldsymbol{\ell}$}] at (0,0) {};

\draw (0,0)--(0.5,1.5)--(-0.3,1)--(0,0);
\draw (0.5,1.5)--(0,2.5)--(1,3);
%\node[shape=circle, fill = black, minimum size = 1.5mm, inner sep = 0pt] at (-0.3,1) {};
%\node[shape=circle, fill = black, minimum size = 1.5mm, inner sep = 0pt] at (0,2.5) {};
\node[shape=circle, fill = black, minimum size = 1.5mm, inner sep = 0pt] at (1,3) {};

\draw (0,0)--(0.333,1)--(0.5,0.5)--(0,0);
\draw (0.333,1)--(0.666,2)--(1,1.5)--(0.333,1);
\draw (0.666,2)--(1,3)--(1.2,2.5)--(0.666,2);
\node[shape=circle, fill = red, minimum size = 1.5mm, inner sep = 0pt] at (0.333,1) {};
\node[shape=circle, fill = red, minimum size = 1.5mm, inner sep = 0pt] at (0.666,2) {};
%\node[shape=circle, fill = black, minimum size = 1.5mm, inner sep = 0pt] at (1,1.5) {};
%\node[shape=circle, fill = black, minimum size = 1.5mm, inner sep = 0pt] at (1.2,2.5) {};
\node[shape=circle, fill = black, minimum size = 1.5mm, inner sep = 0pt] at (1,3) {};
%\node[shape=circle, fill = black, minimum size = 1.5mm, inner sep = 0pt] at (0.5,0.5) {};
\node[shape=circle, fill = red, minimum size = 1.5mm, inner sep = 0pt] at (0.5,1.5) {};

\end{tikzpicture}}
\end{center}
\caption{Pasting algorithm. Step 1: two surfaces $\Gamma_1$ and $\Gamma_2$ are given together with their approximations $\Gamma_{1,h}$ and $\Gamma_{2,h}$. The triangles having an edge on $\ell$ are depicted and their nodes on $\ell$ are black-marked. Step 2: by pasting the triangulated surfaces, a nonconforming polygonation of $\Gamma = \Gamma_1\cup\Gamma_2$ is formed, due to the presence of hanging nodes on $\ell$, that are red-marked.}\label{fig:pastingsurfaces}
\end{figure}

In this section we discuss a possible advantage of SVEM with respect to SFEM. Suppose that $\Gamma$ is made up of two surfaces $\Gamma_1$ and $\Gamma_2$, joining along a straight line, i.e. $\Gamma = \Gamma_1\cup\Gamma_2$ and $\Gamma_1\cap\Gamma_2 = \ell$ is a straight line. Furthermore, suppose we are given two corresponding polygonal surfaces $\Gamma_{1,h}$, $\Gamma_{2,h}$, and that these polygonal surfaces fit $\ell$ exactly, i.e. $\Gamma_{1,h} \cap \Gamma_{2,h} = \ell$. We want to construct a polygonal surface $\Gamma_h$ by pasting $\Gamma_{h,1}$ and $\Gamma_{h,2}$. Such a process is depicted in Fig. \ref{fig:pastingsurfaces} and leads, in general, to a nonconforming overall triangulation. For this reason, the possibility of handling hanging nodes is crucial in this pasting process. It is well-known that the triangular FEMs, including SFEM, are not applicable to nonconforming meshes. For this reason, pasting algorithms for standard FEMs typically need additional deforming and node matching steps, see for instance \cite{kanai1999interactive, sharf2006snappaste}.

\section{Implementation}
\label{sec:implementation}
In this section we will discuss how to implement the SVEM for $k=1$ using only information on the mesh and the nodal values of the load term $f$. We will not consider the case $k\geq 2$ because, as discussed in Section \ref{estimates}, increasing $k$ does not improve the convergence rate of the method.\\
We observe that, from \eqref{zeroaveragedspace} and \eqref{discretezeroaveragedloadterm}, problem \eqref{computablediscreteformulation} is equivalent to
\begin{equation}
\label{implementation_discreteformulation}
\begin{cases}
&u_h\in V_{h}^k\\
&a_h(u_h,\phi_h) = \langle f_h,\phi_h\rangle_h \qquad \forall\ \phi_h\in V_{h}^k\\
&\langle u_h, 1\rangle_{\xLtwo(\Gamma_h)} = 0.
\end{cases}
\end{equation}
Notice that, since $k=1$, the overall number of degrees of freedom is equal to the number of nodal points. Now, for every $i=1,\dots,N$, let $\phi_i\in V_h^1$ be the $i$-th basis function defined by $\text{dof}_j(\phi_i) = \delta_{ij}$, for all $j=1,\dots, N$. We express the numerical solution of \eqref{implementation_discreteformulation} in the basis $\{\phi_i\}_{i=1}^N$ as
\begin{equation*}
u_h(\boldx) = \sum_{j=1}^N \xi_j \phi_j(\boldx),\qquad \forall \boldx\in\Gamma_h,
\end{equation*}
with $\xi_j\in\mathbb{R}$ for all $j=1,\dots,N$. Problem \eqref{computablediscreteformulation} is then equivalent to
\begin{align}
\label{linearsystem1}
&\displaystyle\sum_{j=1}^N a_h(\phi_i,\phi_j) \xi_j = \langle f_h, \phi_i \rangle_h,\qquad \forall i=1,\dots,N,\\
\label{linearsystem2}
&\displaystyle\sum_{j=1}^N \langle 1, \phi_j \rangle_{\xLtwo_h} \xi_j = 0.
\end{align}
Problem \eqref{linearsystem1}-\eqref{linearsystem2} is a rectangular $(N+1)\times N$ linear system that has, from Corollary \ref{crllr:errorestimate}, a unique solution. We want to rephrase this problem as a square $N\times N$ linear system. To this end, notice that the function $\bar{\phi} := \sum_{i=1}^N\phi_i$ fulfills $\text{dof}_j(\bar{\phi}) = 1$ for all $j=1,\dots,N$ and thus, from \eqref{localspace}, we have
\begin{equation}
\label{sumequalto1}
\sum_{i=1}^N\phi_i(\boldx) = 1,\qquad \forall \boldx\in\Gamma_h.
\end{equation}
We show that the sum of all equations in \eqref{linearsystem1} vanishes. In fact, for the left hand side of \eqref{linearsystem1}, using \eqref{consistency} and \eqref{sumequalto1}, we have that
\begin{equation*}
\sum_{i=1}^N\sum_{j=1}^N a_h(\phi_i,\phi_j)  \xi_j= \sum_{j=1}^N  a_h\left(\sum_{i=1}^N\phi_i,\phi_j\right) \xi_j= \sum_{j=1}^N  a_h(1,\phi_j) \xi_j= \sum_{j=1}^N  \bar{a}(1,\phi_j) \xi_j=0,
\end{equation*}
while for the right hand side of \eqref{linearsystem1}, from \eqref{discretezeroaveragedloadterm} and \eqref{sumequalto1} we have
\begin{equation*}
\sum_{i=1}^N\langle f_h, \phi_i \rangle_h = \langle f_h,1\rangle_h = 0.
\end{equation*}
We conclude that the sum of equations \eqref{linearsystem1} vanishes. This implies that we can remove, for instance, the $N$-th equation in \eqref{linearsystem1}. System \eqref{linearsystem1}-\eqref{linearsystem2} is then equivalent to the $N\times N$ system
\begin{equation*}
\begin{cases}
&\displaystyle\sum_{j=1}^N a_h(\phi_i,\phi_j) \xi_j = \langle f_h, \phi_i \rangle_h,\qquad \forall i=1,\dots,N-1,\\
&\displaystyle\sum_{j=1}^N \langle 1,\phi_j \rangle_{\xLtwo_h} \xi_j = 0.
\end{cases}
\end{equation*}
Consider the stiffness matrix $A$, the mass matrix $M$, and the load term $\boldb$ defined by
\begin{align*}
&A = (a_{ij}) := a_h(\phi_i,\phi_j),\qquad \forall i,j=1,\dots,N,\\
&M = (m_{ij}) := \langle \phi_i, \phi_j \rangle_{\xLtwo_h},\qquad \forall i,j=1,\dots,N,\\
&\boldb = (b_i) := \langle f_h, \phi_i \rangle_h,\qquad \forall i=1,\dots,N,
\end{align*}
and, for every $E\in\mathcal{T}_h$, consider the local mass matrix $M_E$ defined by
\begin{equation*}
M_E = (m^E_{ij}) := \langle \phi_i,\phi_j\rangle_{\xLtwo_h,E},\qquad \forall i,j: \boldx_i,\boldx_j\in\text{nodes}(E).
\end{equation*}
The matrices $A$ and $M_E$ and $M$ can be computed as detailed in \cite{beirao2014hitchhiker}. To compute the load vector $\boldb$ we observe that, from \eqref{computablefunctional} and the definition of the basis functions, it holds that
\begin{equation}
\label{implementation_loadterm}
b_i = \sum_{E: \boldx_i\in\text{nodes}(E)} \frac{1}{n_E}\int_E f_h,\qquad \forall i=1,\dots,N.
\end{equation}
We are left to compute the integrals in \eqref{implementation_loadterm} as follows. The nodal values of the load term $f_h$ are computed by
\begin{equation*}
f_h(\boldx_k) = f_I(\boldx_k) - \frac{\langle f_I, 1\rangle_{\xLtwo_h}}{|\Gamma_h|} = f(\boldx_k) - \frac{\sum_{i=1}^N \langle \phi_i, 1\rangle_{\xLtwo_h}f(\boldx_i)}{\langle 1,1\rangle_{\xLtwo_h}} = f(\boldx_k) - \frac{\sum_{i=1}^N (\sum_{j=1}^N m_{ij}) f(\boldx_i)}{\sum_{i=1}^N\sum_{j=1}^N m_{ij}},
\end{equation*}
for all $k=1,\dots,N$. For every $E$, the integral of $f_h$ on $E$ is given by
\begin{equation*}
\int_E f_h \underset{\eqref{exactintegral}}{=} \langle f_h,1 \rangle_{\xLtwo_{h},E} = \sum_{i: \boldx_i\in\text{nodes}(E)} \langle \phi_I,1 \rangle_{\xLtwo_{h},E} f_h(\boldx_i) = \sum_{i: \boldx_i\in\text{nodes}(E)} \left(\sum_{j: \boldx_j\in\text{nodes}(E)}m^{E}_{ij}\right)f_h(\boldx_i).
\end{equation*}
In conclusion, the discretisation of the Laplace-Beltrami equation \eqref{weakform0} by SVEM is given by the following sparse, square, full-rank linear algebraic system
\begin{equation}
\label{finalsystem}
\begin{cases}
&\displaystyle\sum_{j=1}^N a_{ij} \xi_j = b_i,\qquad \forall i=1,\dots,N-1,\\
&\displaystyle\sum_{j=1}^N \left(\sum_{i=1}^N m_{ij}\right) \xi_j = 0.
\end{cases}
\end{equation}
For instance, the matrix form of system \eqref{finalsystem} in \textsc{{\Large M}ATLAB} language is given by
\begin{equation*}
\texttt{xi = [A(1:N-1,:) \hspace*{-2mm}; ones(1,N)*M]}\hspace*{3mm} \verb \ \hspace*{3mm} \texttt{[b(1:N-1); 0].}
\end{equation*}

\section{Numerical examples}
\label{sec:numericalexamples}
In this section we will validate the theoretical findings through numerical experiments.\\
In Experiment 1, a Laplace-Beltrami problem on the unit sphere, approximated with a polygonal mesh, is used as a benchmark to test the convergence rate in \eqref{convergenceestimate}. The experiment also shows the robustness of the method with respect to ``badly shaped'' meshes, i.e with polygons of very different size and very tight, thus confirming the generality of assumptions \ref{itm:firstregularity}-\ref{itm:secondregularity}. In Experiment 2, we show an example of mesh pasting, as discussed in Section \ref{sec:pasting}.
\subsection{Experiment 1}
\begin{figure}
\begin{center}
\subfigure[Polygonal approximation $\Gamma_h$ of the unit sphere $\Gamma$,\newline made up of triangles and hexagons, with meshsize $h= 0.4013$.]{\includegraphics[scale=0.57]{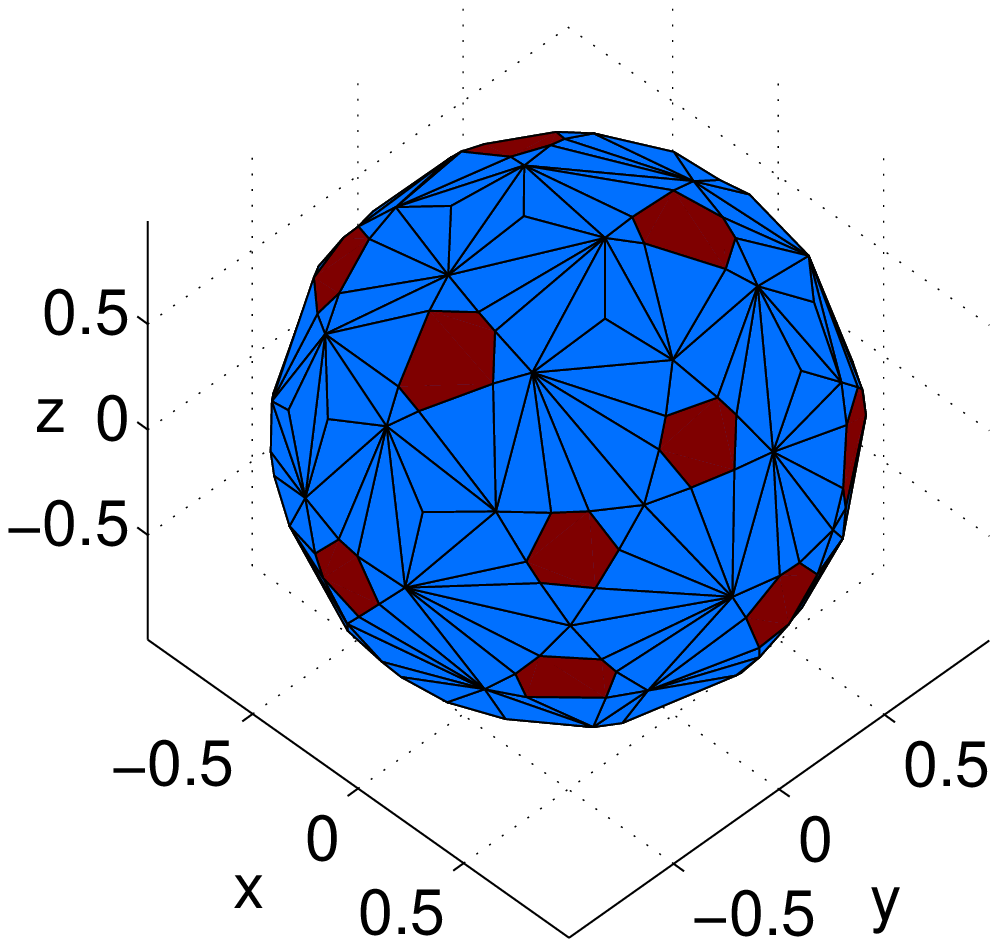}\label{fig:spheremesh}}
\subfigure[Numerical solution obtained on $\Gamma_h$, for $h= 0.0798$.]{\includegraphics[scale=0.57]{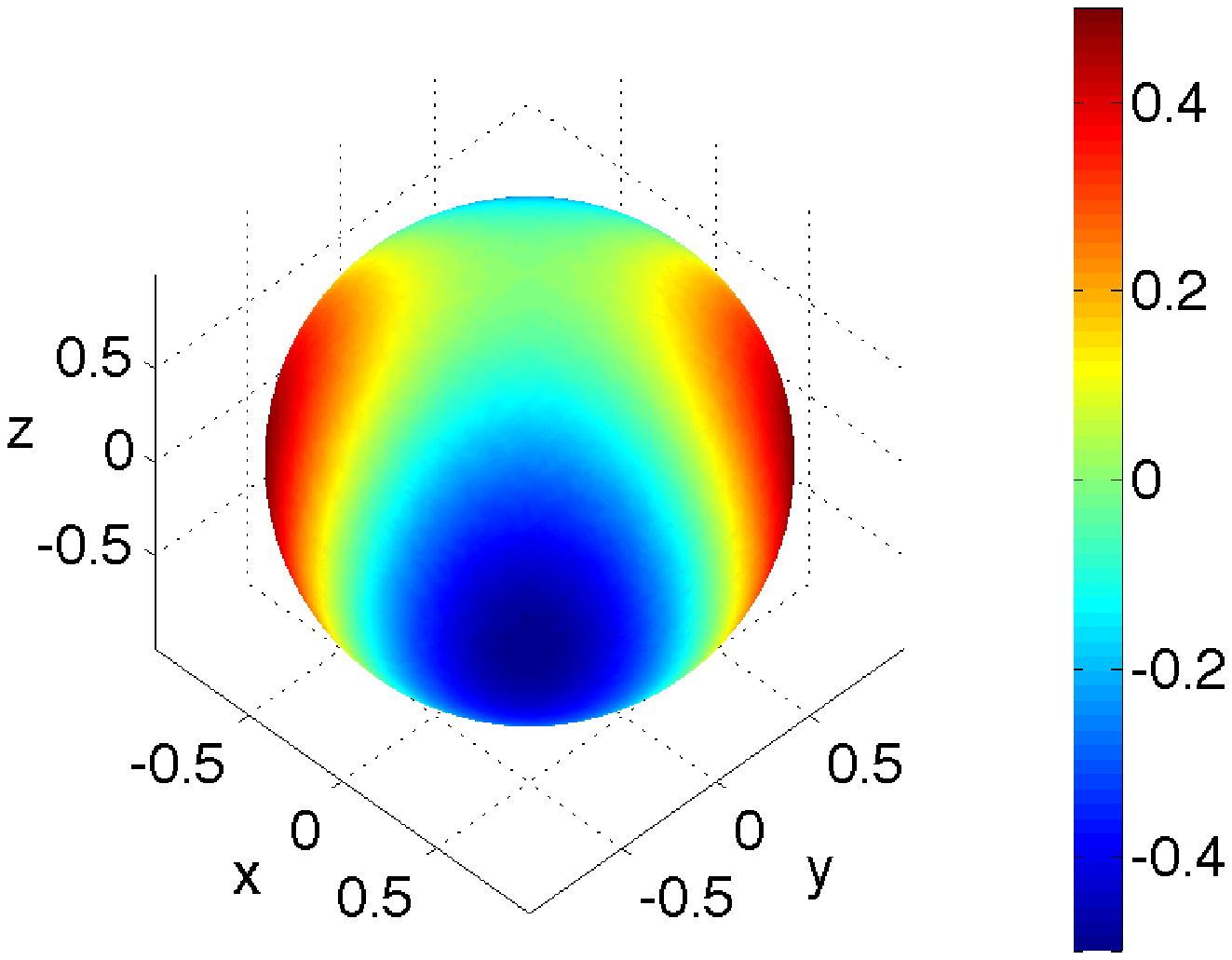}\label{fig:spheresol}}
\caption{(a) The coarsest mesh used for the convergence study; (b) the numerical solution of problem \eqref{experiment1} on the finest mesh.}
\end{center}
\end{figure}

\begin{figure}
\begin{center}
\includegraphics[scale=0.8]{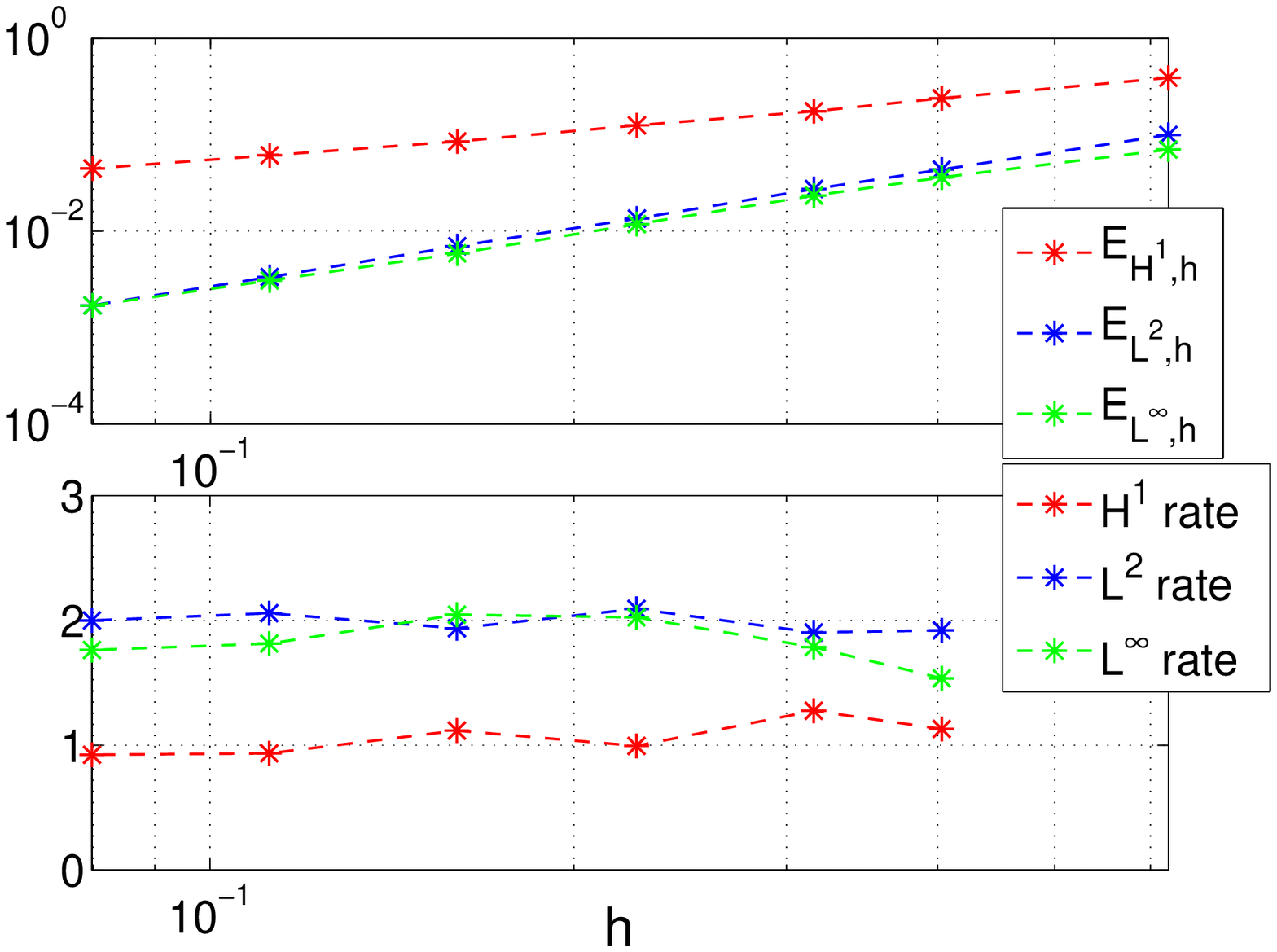}
\caption{Convergence results for problem \eqref{experiment1}.}\label{fig:sphereerror}
\end{center}
\end{figure}

In this experiment we solve the Laplace-Beltrami equation
\begin{equation}
 \label{experiment1}
\begin{cases}
-\hspace{-3mm}&\Delta_\Gamma u(x,y,z) = 6xy, \quad (x,y,z)\in \Gamma,\\
&\int_\Gamma u(x,y,z)\xdif\sigma = 0
\end{cases}
\end{equation}
on the unit sphere $\Gamma := \{ (x,y,z)\in\xR^3 \mid x^2+y^2+z^2=1\}$, whose exact solution is given by
\begin{equation*}
u(x,y,z) = xy,\quad (x,y,z)\in \Gamma.
\end{equation*}
We solve the problem on a sequence of seven polygonal meshes, with decreasing meshsize $h$, made up with triangles and hexagons whose vertices lie on $\Gamma$ and we test the convergence rate as follows. Let $u_I$ be the interpolant, defined in \eqref{interpolationlemmaestimate}, of the exact solution $u$ and let $\delta_h :=u_I-u_h$. We consider the following approximations of the $\xLtwo$, $\xLinfty$ and $\xHone$ errors, respectively:
\begin{align}
\label{approximateL2err}
&\text{E}_{\xLtwo,h} := \langle \delta_h, \delta_h\rangle_{\xLtwo_h};\\
&\text{E}_{\xLinfty,h} := \max_{P\in\text{nodes}(\Gamma_h)}(\delta_h);\\
\label{approximateH1err}
&\text{E}_{\xHone,h} := a_h(\delta_h, \delta_h),
\end{align}
where the forms $a_h(\cdot,\cdot)$ and $\langle \cdot, \cdot\rangle_{\xLtwo_h}$ are defined in \eqref{globalapproximateform} and \eqref{approximateL2form}, respectively. These approximations are $\mathcal{O}(h^2)$-accurate, see for instance \cite{vacca2015virtual}.
The need of defining these approximate norms and seminorms arise from the presence of the virtual basis functions that are not known in closed form. These norms are reminiscent of the approximate $\xLtwo$ norm used for instance in \cite[Equation 46]{vacca2015virtual}, but also account for the fact that, in our case, the exact and the numerical solutions are defined on different surfaces. The convergence rate in the norms and seminorms defined in \eqref{approximateL2err}-\eqref{approximateH1err} is estimated by computing these errors as functions of $h$.\\
The coarsest of the polygonal meshes under consideration (meshsize $h=0.6209$) is shown in Figure \ref{fig:spheremesh}. The numerical solution obtained on the finest mesh (meshsize $h = 0.0798$) is shown in Figure \ref{fig:spheresol}. The convergence results are shown in Fig. \ref{fig:sphereerror}. The convergence is linear in $\xHone$ norm and, even if the method is not designed for optimal $\xLtwo$ and $\xLinfty$ convergence, it appears to be quadratic in $\xLtwo$ norm and almost quadratic in $\xLinfty$ norm. We remark that the considered meshes, like the one in Fig. \ref{fig:spheremesh}, have polygons of very different size and shape, this means that the regularity assumptions \ref{itm:firstregularity} and \ref{itm:secondregularity} are quite weak and the method is thus robust with respect to badly shaped meshes.

\subsection{Experiment 2}
\begin{figure}
\begin{center}
\subfigure[Nonconforming polygonal approximation $\Gamma_h$ of $\Gamma$, with meshsize $h=0.2542$, made up of rectangles and pentagons.]{\includegraphics[scale=0.56]{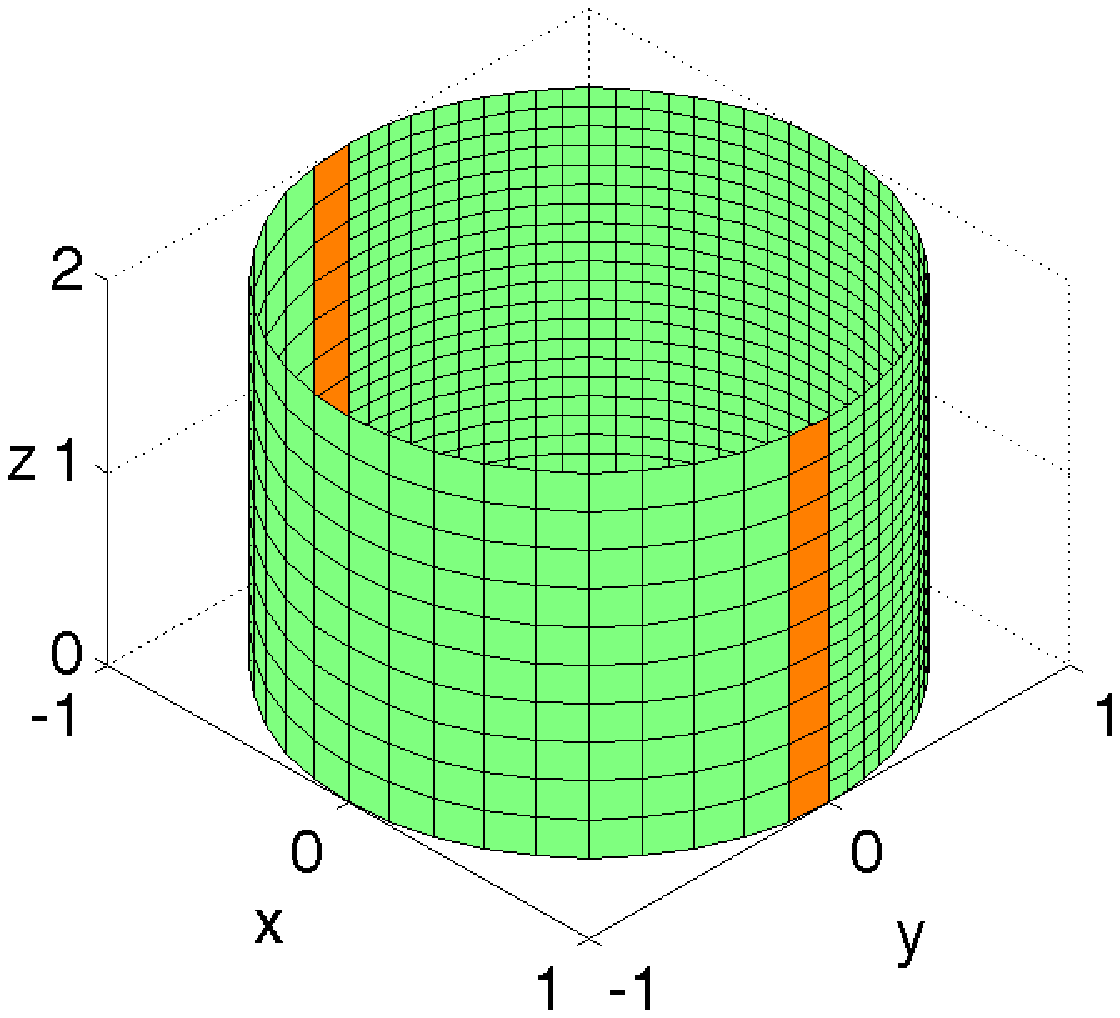}\label{fig:cylindermesh}}
\hspace*{0mm}
\subfigure[Numerical solution obtained on the mesh $\Gamma_h$.]{\includegraphics[scale=0.56]{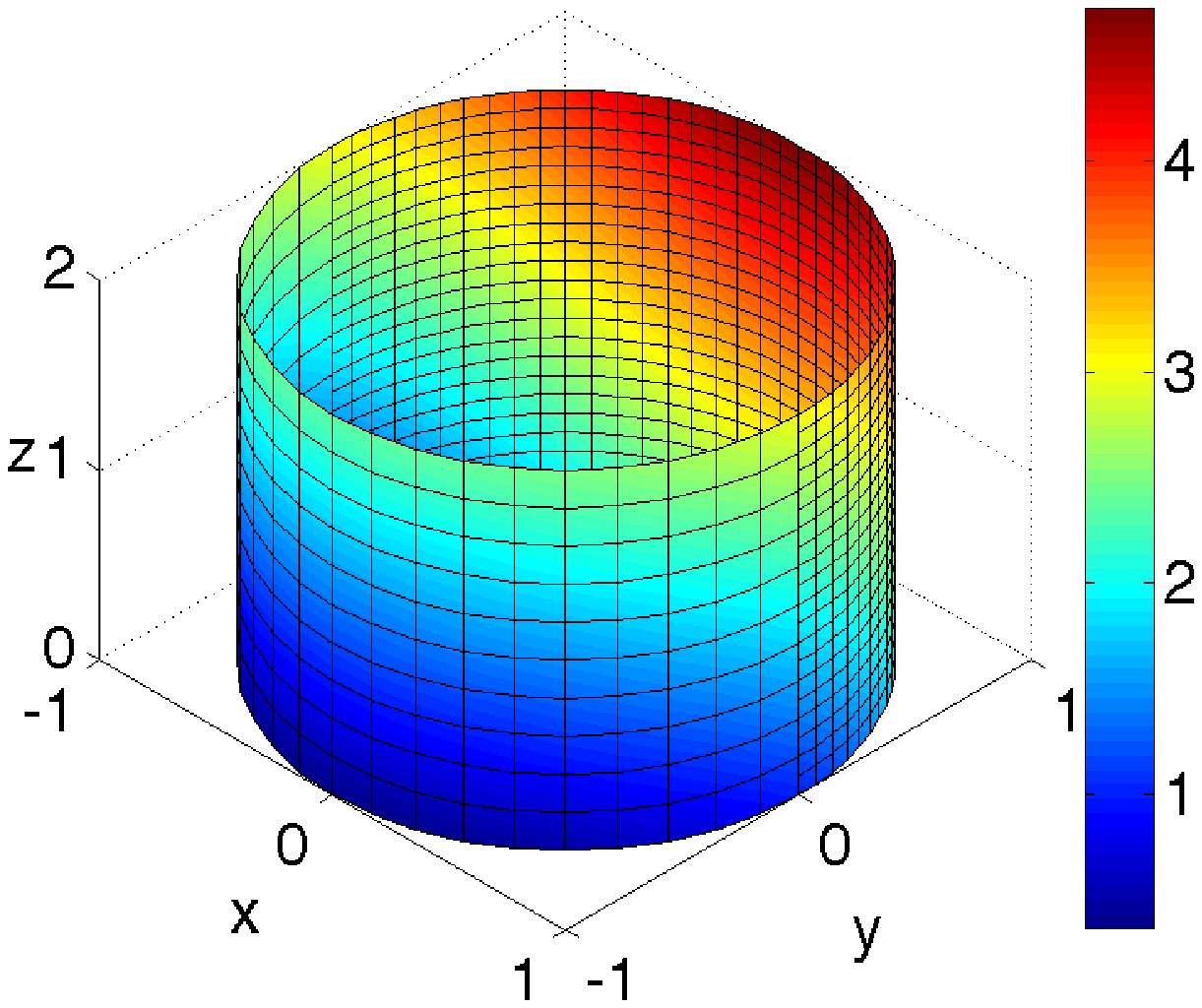}\label{fig:cylindersolution}}
\caption{Polygonal approximation of the cylinder \eqref{cylinder} and corresponding numerical solution of Problem \eqref{pastingproblem}.}
\end{center}
\end{figure}

\begin{figure}
\begin{center}
\includegraphics[scale=0.8]{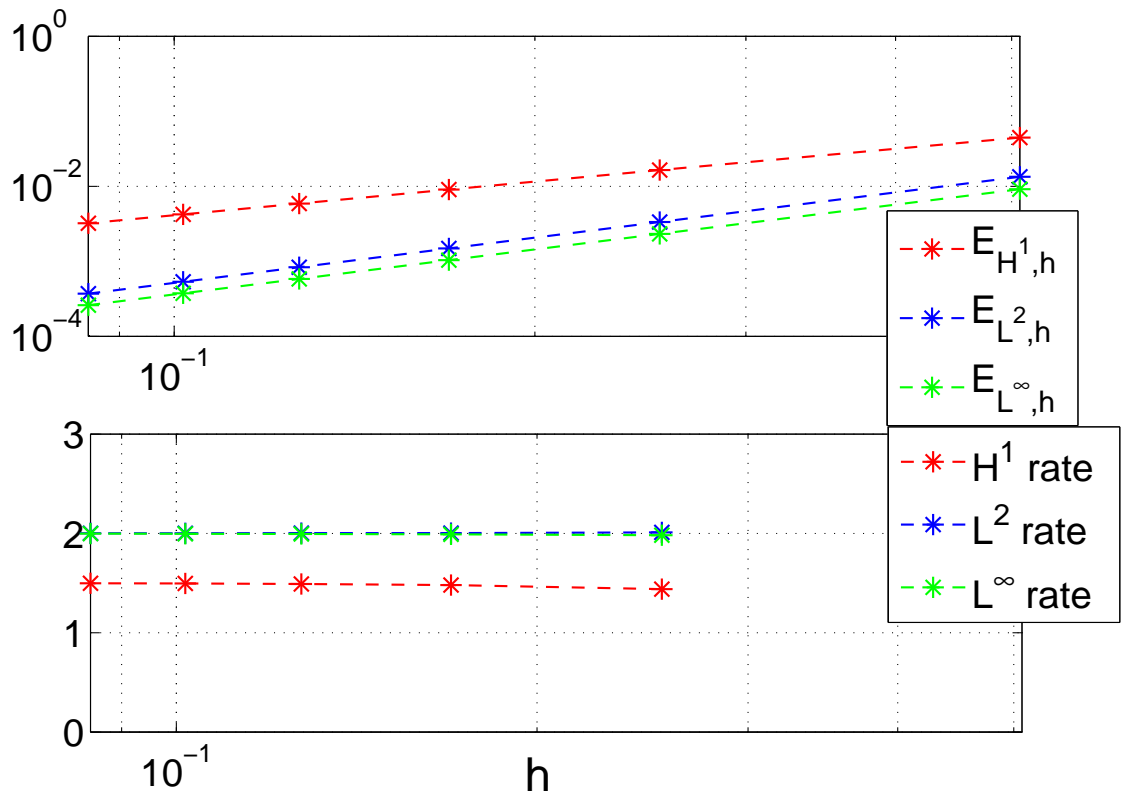}
\caption{Convergence results for problem \eqref{pastingproblem}.}\label{fig:pastingerror}
\end{center}
\end{figure}

In this experiment we solve the Laplace-Beltrami equation and we address the problem of pasting two surfaces. We start by considering the cylinder
\begin{equation}
\label{cylinder}
\Gamma := \{(x,y,z)\in\xR^3 \mid x^2+y^2 =1 \wedge 0\leq z \leq 2\}
\end{equation}
and we split it into two parts
\begin{align*}
\Gamma_1 := \Gamma \cap \{y\geq 0\};\\
\Gamma_2 := \Gamma \cap \{y \leq 0\}.
\end{align*}
We consider the Laplace-Beltrami equation
\begin{equation}
\label{pastingproblem}
\begin{cases}
-\Delta_\Gamma \hspace*{-3mm} &u = (y-x^2)e^y\quad \text{on}\ \Gamma,\\
& u = \bar{u}\quad \text{on}\ \partial \Gamma,
\end{cases}
\end{equation}
where $\bar{u}$ is the exact solution, given by
\begin{equation*}
\bar{u}(x,y,z) = e^y+z,\qquad (x,y,z)\in\Gamma.
\end{equation*}
Notice that the considered surface has a non-empty boundary and the boundary conditions are of homogeneous Dirichlet type, see Remark \ref{rmk:surfaceswithboundary}.
We consider a family of meshes defined as follows. Let $N\in\xN$. The half cylinder $\Gamma_1$ is approximated with $8N^2$ equal rectangular elements constructed on the following $(4N+1)(2N+1)$ gridpoints
\begin{equation*}
P_{ij} = \left(\cos\frac{i}{4N}\pi, \sin\frac{i}{4N}\pi, \frac{j}{N}\right),\qquad i=0,\dots 4N,\quad j= 0,\dots, 2N,
\end{equation*}
while the half cylinder $\Gamma_2$ is approximated with $2N^2$ equal rectangular elements constructed on the following $(2N+1)(N+1)$ gridpoints:
\begin{equation*}
P_{ij} = \left(\cos\left(\frac{i}{2N}+1\right)\pi, \sin\left(\frac{i}{2N}+1\right)\pi, \frac{2j}{N}\right),\qquad i=0,\dots 2N,\quad j= 0,\dots, N.
\end{equation*}
By pasting these meshes we end up with a nonconforming mesh $\Gamma_h$ made up of $10N^2$ elements, of which $2N(5N-1)$ rectangles and $2N$ degenerate pentagons with one hanging node each. For $N=10$, this nonconforming mesh is shown in Fig. \ref{fig:cylindermesh}, in which the rectangles are green and the pentagons are orange while the corresponding numerical solution is shown in Fig. \ref{fig:cylindersolution}.\\
To test the convergence, we consider a sequence of six meshes of the type described above, by increasing $N=5i$, $i=1,\dots,6$. Notice that the relation between $N$ and the meshsize $h$ is
\begin{equation*}
h = 2\sin \left(\frac{\pi}{8N}\right),
\end{equation*}
and thus $h = \mathcal{O}(\frac{1}{N})$. For all $i=1,\dots,6$, the errors in the norms and seminorms \eqref{approximateL2err}-\eqref{approximateH1err} are shown in Fig. \ref{fig:pastingerror} as functions of $h$. The experimental convergence rate is quadratic in the approximate $\xLtwo$ and $\xLinfty$ norms and superlinear in the approximate $\xHone$ seminorm.

\section*{Conclusions}
In this study, we have considered a Surface Virtual Element Method (SVEM) for the numerical approximation of the Laplace-Beltrami equation on smooth surfaces, by generalising the VEM on planar domains \cite{beirao2013basic} and the SFEM \cite{dziuk2013finite}.
By extending the results in \cite{beirao2013basic} and \cite{dziuk2013finite} we have shown, under minimal regularity assumptions for the mesh, optimal asymptotic error estimates (i) for the interpolation in the SVEM function space, (ii) for the approximation of the surface and (iii) for the projection onto the polygonal surface. In particular, the geometric error arising from the approximation of the surface is quadratic in the meshsize and is thus asymptotically dominant with respect to the interpolation error for higher polynomial orders $k>1$. Improving the geometric error is necessary to increase the convergence rate of the overall method. To this end, following \cite{demlow2009higher}, curved polygonal elements could be used. This will be subject of future investigations.
We have also shown the existence and uniqueness of the numerical solution and, for $k=1$, its first order $\xHone$ convergence.
To highlight and advantage of the SVEM technique, we have shown that the process of pasting two meshes along a straight line leads to a nonconforming mesh that can be easily handled by the SVEM. We have presented two numerical examples for the Laplace-Beltrami equation to (i) test the convergence rate of the SVEM method for $k=1$ on the unit sphere and (ii) to show its application to mesh pasting on a cylindrical surface. Since the Laplace-Beltrami equation is endowed with the zero average condition, we have shown, for $k=1$, how to include this condition in the implementation, obtaining a sparse, square, full-rank linear algebraic system, using only information on the mesh and the nodal values of the load term.
The Laplace-Beltrami equation has been considered because it is the simplest PDE defined on a surface. Given the satisfactory behavior of the SVEM on this test problem, the authors believe that it is worth extending the methodology to more complicated surface PDEs. This will be done in future studies.

\end{document}